\documentclass[12pt,a4paper,reqno]{amsart}

\usepackage{graphicx}
\usepackage{latexsym}
\usepackage{amsfonts}

\makeindex

\newcommand{\epsfig}[1]{\reminder{Missing picture here!}}

\newcounter{null}
\newcommand{\Nplus}{{\mathbb N^*}}

\newcommand{\R}{\mathbb {R}}
\newcommand{\C}{\mathbb {C}}

\newfont\bbf{msbm10 at 12pt}
\newcommand\ovl[1]{\overline {#1}}

\newfont\script{eusm10 at 12pt}

\newcommand{\Sym}{\Sigma^1}

\newcommand{\Syms}{\Sigma^\star}

\newtheorem{satz}{\hspace{-1.5mm}}[section]

\newcommand{\lineclear}
      {\rule{0pt}{0pt}\nopagebreak\par\nopagebreak\noindent}
\newenvironment {lemma} [1]
{\begin{satz} {\bf Lemma (#1) }\nopagebreak}{\end{satz}}
\newenvironment {defi} [1]
{\begin{satz} {\bf Definition (#1)}\nopagebreak}{\end{satz}}
\newenvironment {theo} [1]
{\begin{satz} {\bf Theorem (#1)}\nopagebreak}{\end{satz}}
\newenvironment {prop} [1]
{\begin{satz} {\bf Proposition (#1) }\nopagebreak}{\end{satz}}
\newenvironment {coro} [1]
{\begin{satz} {\bf Corollary (#1) }\nopagebreak}{\end{satz}}
\newenvironment {coro*}
{\begin{satz} {\bf Corollary }\nopagebreak}{\end{satz}}
\newenvironment {example} [1]
{\begin{satz} {\bf Example (#1)}\rm\nopagebreak }{\end{satz}}

\renewcommand{\proof}{\medskip\par\noindent{\sc Proof. }}
\newcommand{\proofof}[1]{\medskip\par\noindent{\sc Proof of #1}.}

\newcommand\nix{\rule{0pt}{2pt}}
\renewcommand{\qed}{\qedd\par\medskip\noindent}
\newcommand{\qedd}{\nix\nolinebreak\hfill\hfill\nolinebreak$\Box$}

\newcommand{\remark}{\par\medskip \noindent {\sc Remark. }}
\newcommand{\Intro}[1] {{\small #1}\par\medskip}

\newcommand{\reminder}[1]{{\sf #1} }
\newcommand{\hide}[1]{}

\newcommand{\evil}{evil}
\newcommand{\tame}{tame}
\newcommand{\Evil}{Evil}
\newcommand{\Tame}{Tame}

\newcommand{\0}{{\tt 0}}
\newcommand{\1}{{\tt 1}}
\renewcommand{\*}{{\tt \star}}
\newcommand\A{{\mathcal A}}
\newcommand\Abar{\ovl{\A}}

\newcommand{\IntAdr}{\to}

\renewcommand{\text}[1]{\mbox{#1}}  
\newcommand{\ph}{\varphi}
\renewcommand\phi{\varphi}

\renewcommand\theta{\vartheta}

\newcommand{\z}{\zeta}

\newcommand{\orb}{\mbox{\rm orb}}

\newcommand\sm{\setminus}

\newcommand{\step}{\mbox{{\sc Step}}}  

\newcounter{stepcount}
    {\begin{list}{Step~\arabic{stepcount}.}{\usecounter{stepcount}}}%
    {\end{list}}

\def\reminder#1{{\sf #1} }
\def\hide#1{}

\hide{
{\begin{satz} {\bf Theorem (\boldmath#1\unboldmath)}
\nopagebreak}{\end{satz}}
}

\begin{document}
\title[Kneading sequences and structure of the Hubbard tree. \today]
{Admissibility of kneading sequences and structure\\[2mm] of Hubbard trees for quadratic polynomials}

\author{Henk Bruin}
\author{Dierk Schleicher}
\address{
Department of Mathematics, 
University of Surrey, 
Guildford GU2 7XH, 
United Kingdom}
\email{H.Bruin@surrey.ac.uk}

\address{
School of Engineering and Science,
Jacobs University Bremen,
P.O. Box 750 561, D-28725 Bremen,
Germany}
\email{dierk@jacobs-university.de}
\hide{
\texttt{http://www.jacobs-university.de/directory/dschleicher/} 
}


\begin{abstract}
Hubbard trees are invariant trees
connecting the points of the critical orbits of postcritically finite
polynomials. Douady and Hubbard \cite{Orsay} introduced these trees
and showed that they encode the essential information of Julia sets in a combinatorial way. The itinerary of the critical orbit within the Hubbard tree is encoded by a (pre)periodic sequence on $\{\0,\1\}$ called \emph{kneading sequence}.

We prove that the kneading sequence completely encodes the Hubbard
tree and its dynamics, and we show how to reconstruct the tree and in
particular its branch points (together with their periods, their
relative posititions, their number of arms and their local dynamics)
in terms of the kneading sequence alone.

Every kneading sequence gives rise to an abstract Hubbard tree,
but not every kneading sequence occurs in real dynamics or in complex
dynamics. Milnor and Thurston~\cite{MT} classified which kneading
sequences occur in real dynamics; we do the same for complex dynamics
in terms of a complex \emph{admissibility condition}.Ê
This complex admissibility condition fails if and only if the abstract Hubbard tree
has a so-called \emph{evil} periodic branch point that is incompatible with
local homeomorphic dynamics on the plane.
\end{abstract}

\maketitle

\section{Introduction}

In complex dynamics, a frequent observation is that many dynamical properties can be encoded in symbolic terms. Douady and Hubbard \cite{Orsay} discovered that Julia sets of polynomial Julia sets could completely be described in terms of a tree that is now called the \emph{Hubbard tree} (at least in the case of postcritically finite polynomials; a complete classification was later given in \cite{BFH,Poirier}). 

We investigate Hubbard trees of postcritically finite quadratic polynomials. We show that these trees can completely be described by a single periodic binary sequence called \emph{kneading sequence} which encodes the location of the critical orbit within the tree. More precisely, we show that all endpoints and all branch points of a Hubbard tree are completely encoded by the kneading sequence, and that these suffice to describe the Hubbard tree and its dynamics up to a natural equivalence relation. We show that orbits of branch points come in two kinds which we call \emph{tame} and \emph{evil}.

Kneading sequences are ubiquitous in real and complex dynamics and they have been studied by many people. Milnor and Thurston \cite{MT} classified all kneading sequences that arise in real dynamics, especially by real quadratic polynomials. We answer the corresponding complex question and classify all kneading sequences that arise in complex dynamics: our admissibility condition
is given in Definition~\ref{DefAdmissCond}. For this, it suffices to restrict attention to sequences that we call \emph{$\*$-periodic}: it turns out every $\*$-periodic kneading sequence is realized by an essentially unique abstract Hubbard tree; so in order to determine which kneading sequences are realized, we can investigate the associated abstract Hubbard trees. 
We point out that our trees are abstract in the sense that they do not come
with an embedding into the complex plane, in contrast to the original definition
of Douady and Hubbard.

Both real and complex admissibility of abstract Hubbard trees are encoded in their branch points: an abstract Hubbard tree is real admissible if it has no branch point at all (the tree is an interval and can be embedded into $\R$); an abstract Hubbard tree is complex admissible if it can be embedded into $\C$ so that the embedding respects the circular order. In terms of our classification of branch points, this means that all branch points of the tree must be tame: every evil branch point is an obstruction to complex admissibility of a kneading sequence, and evil branch points are the only possible obstructions. 
Readers familiar with Thurston's classification \cite{DHT}
of rational maps may see similarities with
obstructions in that classification. In both cases, a combinatorial obstruction 
prevents a branched cover, or a Hubbard tree, from being realized by a 
holomorphic map, in particular by a quadratic polynomial.
In fact, our results are closely connected to Thurston's theorem (even though we do not use it). 

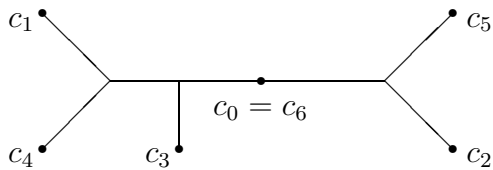
\begin{figure}[htb]
\begin{center}
\begin{minipage}{90mm}
\unitlength=9mm
\begin{picture}(10,3.2) \let\ts\textstyle
\put(3,2){\line(1,0){4}} \put(2,1){\line(1,1){1}}
\put(2,3){\line(1,-1){1}} \put(4,2){\line(0,-1){1}}
\put(7,2){\line(1,1){1}} \put(7,2){\line(1,-1){1}}

\put(5.2,2){\circle*{0.1}} \put(4.5,1.5){$c_0 = c_6$}
\put(2,3){\circle*{0.1}} \put(1.5,2.8){$c_1$}
\put(8,1){\circle*{0.1}} \put(8.2,0.8){$c_2$}
\put(4,1){\circle*{0.1}} \put(3.5,0.8){$c_3$}
\put(2,1){\circle*{0.1}} \put(1.5,0.8){$c_4$}
\put(8,3){\circle*{0.1}} \put(8.2,2.8){$c_5$}
\end{picture}
\end{minipage}
\vskip-20pt
\end{center}
\caption{The Hubbard tree for $1 \IntAdr 2 \IntAdr 4 \IntAdr 5
\IntAdr 6$ contains an {\evil} orbit of period $3$.}
\label{FigEvil}
\end{figure}

The simplest example of a non-admissible sequence is
$\nu = \ovl{\1\0\1\1\0\*}$.
Here the Hubbard tree, shown in Figure~\ref{FigEvil}, has a period $3$ branch point, but the third iterate of $f\colon T \to T$ fixes one arm and permutes the other two transitively. Such branch points cannot be embedded into the plane so that the dynamics respects the circular order of the arms: this is an example of an evil branch point. 

The admissibility condition also applies to kneading sequences that are not preperiodic or $\*$-periodic.
For the interpretation using evil orbits, finite Hubbard trees have to be replaced by dendrites (such as those constructed by Penrose~\cite{Pe}); another interpetation of complex admissibility is in terms of whether the kneading sequence is realized by angle doubling on the circle. Details, and many further properties of Hubbard trees, are the subject of a forthcoming monograph \cite{BKS}.

While Hubbard trees are very good for describing individual Julia sets, it is not quite so easy to tell which trees are close to each other so as to obtain a topology on the space of Hubbard trees. Kneading sequences are helpful here: the natural topology on the space of kneading sequences describes dynamical proximity of Hubbard trees in a way that is compatible, for example, with their location within the Mandelbrot set \cite[Section~6]{BKS}.

Kneading sequences can be recoded in ``human-readable form'' in the form of internal addresses (see Definition~\ref{DefRho} below): in this form, they allow to read off the location in parameter space of any quadratic polynomial just in terms of the kneading sequence \cite{IntAdr,IntAdrNew}, and they help to establish fundamental properties of the Mandelbrot set \cite{ExtRayMandel,Spiders}.

Since all trees in this paper are abstract Hubbard trees, we omit the word ``abstract'' from now on; one should keep in mind that our definition differs from that by Douady and Hubbard in the fact that their trees always come with an embedding into $\C$. Some of our trees cannot be embedded into the plane in a way that is compatible with the dynamics (those which have evil orbits), while others may have many essentially different such embeddings: such trees are realized by several quadratic polynomials with topologically conjugate dynamics.

The structure of the paper is as follows. In Section~\ref{SecTrees}, we define Hubbard trees and fundamental concepts from symbolic dynamics, including itineraries and kneading sequences. In the rest of the paper, we investigate the Hubbard tree associated to a given $\*$-periodic kneading sequence. Existence and uniqueness of this Hubbard tree are shown in \cite{TreeExistence}. We do not assume these results here, so the present paper is essentially self-contained: we investigate properties of trees that we assume to exist (but knowing the existence of all these trees reassures us that we are not investigating empty sets). Section~\ref{SecPeriodic} contains an investigation of all branch points in Hubbard trees as well as the definition of tame and evil branch points and the proof that the embedding of a tree into the plane respecting the dynamics is possible if and only if all periodic orbits are tame. The final section~\ref{SecAdmissCondition} shows how to determine branch points, their number of arms and their type (evil or tame) in terms of the kneading sequence. We also give a constructive uniqueness proof of Hubbard trees which follows from our investigation of the trees in Corollary~\ref{CorUniqueness}.

{\bf Acknowledgement.} 
We gratefully acknowledge that this research was partially supported by the European Marie Curie Research Training Network CODY and by the ESF Research Networking Programme HCAA.  HB also acknowledges support by EPSRC grant GR/S91147/01.

\section{Hubbard Trees}
\label{SecTrees}

In this section, we define Hubbard trees as abstract trees with
dynamics and show their most fundamental properties. Our trees do
not necessarily come with an embedding into the complex plane.

\begin{defi}{Trees, Arms, Branch Points and Endpoints}
\label{DefTree} \lineclear
A {\em tree}\index{tree} $T$ is a finite connected graph without loops.
For a point $x \in T$,  the {\em (global)
arms}\index{arm}\index{global arm}
of $x$ are the connected components of $T\setminus\{x\}$. A
{\em local arm}\index{local arm}
at $x$ is an intersection of a global arm with a
sufficiently small neighborhood of $x$ in $T$. The point $x$ is an
{\em endpoint}\index{endpoint} of $T$ if it has only one arm;
it is a {\em branch point}\index{branch point} if it has at least three arms.
\end{defi}

\noindent
Between any two points $x,y$ in a tree, there exists a unique
closed arc connecting $x$ and $y$; we denote it by $[x,y]$ and its
interior by $(x,y)$.

\begin{defi}{The Hubbard Tree}
\label{DefHubbard} \lineclear
A {\em Hubbard tree}\index{Hubbard tree} is a tree $T$ equipped with a map
$f\colon T\to T$ and a distinguished point, the {\em critical
point}\index{critical point}, satisfying the following conditions:
\begin{enumerate}
\item
$f:T \to T$ is continuous and surjective;
\item
every point in $T$ has at most two inverse images under $f$;
\item
at every point other than the critical point, the map $f$ is a
local homeomorphism onto its image;
\item
all endpoints of $T$ are on the critical orbit;
\item
the critical point is periodic or preperiodic, but not fixed;
\item
(expansivity)\index{expansivity}
if $x$ and $y$ with $x \neq y$ are branch points or
points on the critical orbit, then there is an $n \geq 0$ such
that $f^{\circ n}([x, y])$ contains the critical point.
\end{enumerate}
\end{defi}

We denote the critical point by $c_0=0$ and its orbit by
$\orb_f(c_0) = \{ 0, c_1, c_2, \dots \}$. The {\em critical value}
$c_1$ is the image of the critical point.
We use a standing assumption that $c_1\neq c_0$ in order to avoid having
to deal with counterexamples when the entire tree is a single point.
The branch points and the
points on the critical orbit (starting with $c_0$) will be called
{\em marked points}\index{marked point}.
Notice that the set of marked points is finite
and forward invariant because the number of arms at any point can
decrease under $f$ only at the critical point.

Two Hubbard trees $(T,f)$ and $(T',f')$ are {\em equivalent} if
there is a bijection between their marked points  which is
respected by the dynamics, and if the edges of the tree connect the
same marked points. This is weaker than a topological conjugation.
In particular, we do not care about details of the dynamics between
marked points; there may be intervals of periodic points,
attracting periodic points, and so on. (This is related to an equivalence class
of branched covers in the sense of Thurston as in \cite{DHT,Spiders}.)

\begin{lemma}{The Hubbard Tree}
\label{LemTree} \lineclear
The critical value $c_1$ is an endpoint, and
the critical point $c_0$ divides the tree into at most two parts.
Each branch point is periodic or preperiodic,  it never maps
onto the critical point, and the number of arms is constant along
the periodic part of its orbit. Any arc which does not contain the
critical point in its interior maps homeomorphically onto its
image.
\end{lemma}

\proof
Suppose that $c_1$ has at least two arms. The points $c_2$, $c_3$,
\ldots \, also have at least two arms as long as $f$ is a local
homeomorphism near this orbit. If this is no longer the case at some point,
then the orbit has reached the critical point, and the next image is $c_1$ again.
In any case, all points on the critical orbit have
at least two arms. This contradicts the assumption that all
endpoints of a Hubbard tree are on the critical orbit. Hence $c_1$ has exactly one
arm, and $c_0$ has at most two arms (or its image would not be an
endpoint).

Since near every non-critical point, the dynamics is a local
homeomorphism onto the image, every branch point maps onto a branch
point with at least as many arms. Since the critical point has at
most two arms, it can never be the image of a branch point. The
tree and thus the number of branch points is finite, so every branch
point is preperiodic or periodic and its entire orbit consists of
branch points; the number of arms is constant along the periodic
part of the orbit.

Let $\gamma$ be an arc within the tree. Since $f$ cannot be
constant on $\gamma$ and there is no loop in the tree, the subtree
$f(\gamma)$ has at least two endpoints. If an endpoint of
$f(\gamma)$ is not the image of an endpoint of $\gamma$, then it
must be the image of the critical point since $f$ is a local
homeomorphism elsewhere, and the critical point $0$ must be 
in the interior of $\gamma$.
\qed

In a Hubbard tree $T$ with critical point $c_0$, the set $T\sm\{c_0\}$ consists of at most two connected components; let $T_\1$ be the component containing the critical value and $T_\0=T\sm(T_\1\cup\{c_0\})$ (the set $T_\0$ may or may not be empty). Writing $T_\*=\{c_0\}$, we can define itineraries in the usual way as sequences over $\{ \0, \*, \1\}$. The itinerary $\nu=\nu_1\nu_2\nu_3\dots$ of the critical value $c_1$ is called the \emph{kneading sequence}; it always starts with $\1$.  If $c_0$ is periodic, say of period $n$, then $\nu_n = \*$ and $\nu = \ovl{\1\nu_2\dots\nu_{n-1}\*}$; we call such sequences \emph{$\*$-periodic}.

Write $\Nplus = \{ 1,2,3,\dots\}$ and let
$\Syms$ be the set of all $\nu \in \{ \0, \1\}^{\Nplus}$ and all $\*$-periodic sequences, always subject to the condition that all sequences start with $\nu_1 = \1$.

\begin{defi}{$\rho$-Function and Internal Address}\index{$\rho$-function}
\label{DefRho} \lineclear
For a sequence $\nu \in\Syms$, define
\[
\rho_{\nu}:\Nplus \to \Nplus\cup\{\infty\}, \quad
\rho_\nu(n) = \inf \{ k > n: \nu_k \neq \nu_{k-n} \}.
\]
We usually write $\rho$ for $\rho_\nu$.
For $k\geq 1$, we call
\[
\orb_\rho(k):= k\to\rho(k)\to\rho^{\circ 2}(k)\to\rho^{\circ 3}(k)
\to\dots
\]
the $\rho$-orbit of $k$. The case $k=1$ is the most important one;
we call
\[
\orb_\rho(1) =
1 \to \rho(1) \to \rho^{\circ 2}(1) \to \rho^{\circ 3}(1) \to \dots
\]
the {\em internal address}\index{internal address} of $\nu$.
For real unimodal maps, the numbers $\rho^{\circ k}(1)$
are known as the \emph{cutting times} of the map.
If $\rho^{\circ k+1}(1) = \infty$, then we say that the internal
address is finite: $1 \to \rho(1) \to \ldots \to \rho^{\circ k}(1)$;
as a result, the orbit $\orb_\rho$ is a finite or infinite sequence
that never contains $\infty$.
\end{defi}

The following combinatorial lemma will be used to locate
the images of certain closest precritical points
in Hubbard trees. The proof can be found in \cite{TreeExistence, BKS},
and, with entirely different terminology, in the thesis of Penrose
\cite[Theorem 4.5.3 and Corollary 2.5.3.1]{Pe}.

\begin{lemma}{Combinatorics of $\rho$-Orbits}
\label{LemExactPeriod}\lineclear
Let $\nu\in\Sym$ (not containing a $\*$) and let $m$ belong to the
internal address of $\nu$.
\begin{enumerate}
\item
If $s$ is such that $s < m < \rho(s)$,
then $\orb_{\rho}( \rho(m-s) - (m-s) ) \owns m$.
\hide{  
\item
If $\rho(m) > 2m$, then for every $s \in\{1, \dots, m\}$,  either $m$ or $2m$
belongs to $\orb_{\rho}(s)$.
}
\item
If $\rho(m) = \infty$, then $m$ is the exact period
of $\nu$.
\end{enumerate}
\end{lemma}

\section{Periodic Orbits on Hubbard Trees}
\label{SecPeriodic}


\Intro{ In this section, we discuss periodic points of Hubbard
trees, in particular branch points, and show that they come in two
kinds: {\em {\tame}} and {\em {\evil}}. This determines whether or not
Hubbard trees and kneading sequences are admissible: they are
if and only if there is no {\evil} orbit.}

\begin{lemma}{Characteristic Point}
\label{LemCharacteristic} \lineclear
Let $(T,f)$ be the Hubbard tree with kneading sequence $\nu$.
Let $\{z_1,z_2,\ldots,z_{n}=z_0\}$ be a periodic orbit
which contains no endpoint of $T$. If the critical orbit is
preperiodic, assume also that the itineraries of all points $z_k$
are different from the itineraries of all endpoints of $T$.

Then there are a unique point $z\in\{z_k\}_{k=1}^n$ and two
different components of $T\sm\{z\}$ such that the critical value is
contained in one component and $0$ and all other points $z_k\neq
z$ are in the other one.
\end{lemma}

\begin{defi}{Characteristic Point}
\label{DefCharacteristic} \lineclear
The point $z$ in the previous lemma is called the {\em
characteristic point}\index{characteristic point}
of the orbit $\{z_k\}$; we will always
relabel the orbit cyclically so that the characteristic point is
$z_1$.
\end{defi}

\proofof{Lemma~\ref{LemCharacteristic}}
Note first that every $z_k\neq 0$ (or $z_{k+1}=c_1$ would be an
endpoint).
For each $z_k$, let $X_k$ be the union of all components of
$T\setminus\{ z_k\}$ which do not contain the critical point.
Clearly $X_k$ is non-empty and $f|_{X_k}$ is injective. If $X_k$
contains no immediate preimage of $0$, then $f$ maps $X_k$
homeomorphically into $X_{k+1}$. Obviously, if $X_k$ and $X_l$
intersect, then either $X_k \subset X_l$ or $X_l \subset X_k$.
At least one set $X_k$ must contain an
immediate preimage of $0$:
if the critical orbit is periodic, then every endpoint of $T$
eventually iterates onto $0$, and every $X_k$ contains an endpoint.
If the critical orbit is preperiodic, we need the extra hypothesis
on the itinerary of the orbit $(z_k)$:
if no $X_k$ contains a point which ever iterates to $0$,
then all endpoints of $X_k$ have the same itinerary as $z_k$ in
contradiction to our assumption.

If $X_k$ contains an immediate preimage $w$ of $0$, then the
corresponding $z_k$ separates $w$ from the critical point, i.e., $z_k\in[w,0]$
and thus $z_{k+1}\in[0,c_1]$ (always
taking indices modulo $n$), hence $c_1\in X_{k+1}$.

Among the non-empty set of points $z_{k+1}\in[0,c_1]$, there is a
unique one closest to $c_1$; relabel the orbit
cyclically so that this point is $z_1$. We will show that this is
the characteristic point of its orbit.

For every $k$, let $n_k$ be the number of points from
$\{ z_i \}$ in $X_k$. If $X_k$
does not contain an immediate preimage of $0$, then $n_{k+1} \geq
n_k$. Otherwise, $n_{k+1}$ can be smaller than $n_k$, but
$z_{k+1}\in[0,c_1]$; since no $z_k\in (z_1,c_1]$, we have
$z_{k+1}\in[0,z_1]$ and either $z_{k+1}=z_1$ or $n_{k+1}\geq 1$.

Therefore, if $n_1\geq 1$, then all $n_k\geq 1$; however,
the nesting property of the $X_k$ implies that there is at least
one `smallest' $X_k$ which contains no further $X_{k'}$ and thus no
$z_{k'}$; it has $n_k=0$. Therefore, $n_1=0$; this means that all
$z_k\neq z_1$ are in the same component of $T\sm\{z_1\}$ as $0$.
Since $c_1\in X_1$, the point $z_1$ is characteristic.
\qed

\begin{prop}{Images of Global Arms}
\label{PropBadBranchPoint} \lineclear
Let $z_1$ be a characteristic periodic point of exact period $m$
and let $G$ be a global arm at $z_1$. Then either $0 \notin f^{\circ k}(G)$
for $0 \leq k < m$ (and in particular the first return map of $z_1$ maps $G$
homeomorphically onto its image), or the first return map of $z_1$ sends
the local arm in $G$ to the local arm at $z_1$ pointing to the critical point or
the critical value.
\end{prop}
\proof
Let $z_k:=f^{\circ(k-1)}(z_1)$ for $k\geq 1$.
Consider the images $f(G)$, $f(f(G))$, etc.\ of the global arm $G$;
if none of them contains $0$ before $z_1$ returns to
itself, then $G$ maps homeomorphically onto its image under the
first return map of $z_1$ and the claim follows. Otherwise, there is a first index $k$
such that $f^{\circ(k-1)}(G) \ni 0$, so that the image arm at $z_k$
points to $0$; so far, the map is homeomorphic on $G$.
If $z_k=z_0$, then the image point is $z_1$ and
the local image arm at $z_1$ points to $c_1$. If $z_k\neq
z_0$, then the local arm at $z_k$ points to $0$ and the image
arm at $z_{k+1}$ points to $c_1$; since $z_1$ is
characteristic, the image arm points also to $z_1$. Continuing the
iteration, the image arms at the image points will always
point to some $z_l$.
When the orbit finally reaches $z_0$, the local arm points to some
$z_{l'}$. If it also points to $0$, then the image at $z_1$ will
point to $c_1$ as above; otherwise, it maps homeomorphically and
the image arm at $z_1$ points to $z_{l'+1}$.
By Lemma~\ref{LemCharacteristic}, the only
such arm is the arm to the critical point.
\qed

\begin{coro}{Two Kinds of Periodic Orbits}
\label{CorTwoKindsOrbits} \lineclear
Let $z_1$ be the characteristic point of a periodic orbit of branch
points. Then the first return map either permutes all the local
arms transitively, or it fixes the arm to $0$ and permutes
all the other local arms transitively.
\end{coro}
\proof
Let $n$ be the exact period of $z_1$. Since the periodic orbit does
not contain the critical point by Lemma~\ref{LemTree}, the map
$f^{\circ n}$ permutes the local arms of $z_1$. Let $G$ be any
global arm at $z_1$. It must eventually map onto the critical
point, or the marked point $z_1$ would have the same itinerary as
all the marked points in $G$, contradicting the expansivity
condition. By Lemma~\ref{PropBadBranchPoint}, the orbit of any
local arm at $z_1$ must include the arm at $z_1$ to $0$ or to $c_1$
or both, and there can be at most two orbits of local arms.

Consider the local arm at $z_1$ to $0$. The corresponding global
arm cannot map homeomorphically, so $f^{\circ n}$ sends this local
arm to the arm pointing to $0$ or to $c_1$. If the image local arm
points to $c_1$, then all local arms at $z_1$ are on the same
orbit, so $f^{\circ n}$ permutes these arms transitively.
If $f^{\circ n}$ fixes the local arm at $z_1$ pointing
to $0$, then the orbit of every other local arm must include the arm
to $c_1$, so all the other local arms are permuted transitively.
\qed

\begin{defi}{{\Tame} and {\Evil} Orbits}
\label{DefTameEvil}\lineclear
A periodic orbit of branch points is called {\em {\tame}}\index{{\tame}}
if all its local arms are on the same cycle, and it is called {\em
{\evil}}\index{{\evil}} otherwise.
\end{defi}
\remark
Obviously, {\evil} orbits are characterized  by the property that not
all local arms have equal periods; their first return dynamics is
described in Corollary~\ref{CorTwoKindsOrbits}. For periodic points
(not containing a critical point)
with two local arms, the situation is analogous: the first return
map can either interchange these arms or fix them both. It will
become clear below that periodic points with only two arms are less
interesting than branch points; however,
Proposition~\ref{PropBranchType} shows that they have similar
combinatorial properties.
\hide{
In Section~\ref{SecParameter}, we will
also call a periodic point with two arms {\em {\tame}} if the first
return map permutes the two local arms (but we reserve the
adjective {\em {\evil}} for branch points because only the latter
destroy admissibility; see Proposition~\ref{PropEmbedding}).
}

\begin{lemma}{Global Arms at Branch Points Map Homeomorphically}
\label{LemHomeo} \lineclear
Let $z_1$ be the characteristic point of a periodic orbit of period
$n$ and let $q\geq 3$ be the number of arms at each point. Then the
global arms at $z_1$ can be labelled $G_0$, $G_1$, \ldots, $G_{q-1}$
so that $G_0$ contains the critical point, $G_1$ contains the
critical value, and the arms map as follows:
\begin{itemize}
\item
if the orbit of $z_1$ is {\tame}, then the local arm
$L_0 \subset G_0$ is
mapped to the local arm $L_1 \subset G_1$ under
$f^{\circ n}$; the global arms
$G_1,\ldots, G_{q-2}$ are mapped homeomorphically onto their images
in $G_2,\ldots, G_{q-1}$, respectively, and the local arm
$L_{q-1}\subset G_{q-1}$ is mapped to $L_0$;
\item
if the orbit is {\evil}, then the local arm $L_0$ is fixed under
$f^{\circ n}$, the global arms $G_1,\ldots, G_{q-2}$ are mapped
homeomorphically onto their images in $G_2,\ldots, G_{q-1}$,
respectively, and the local arm $L_{q-1}\subset G_{q-1}$ is sent to
the local arm $L_1$; however, the global arm $G_{q-1}$ maps onto
the critical point before reaching $G_1$.
\end{itemize}
In particular, if the critical orbit is periodic, then its period must strictly exceed
the period of any periodic branch point.
\end{lemma}

\proof
We will use Proposition~\ref{PropBadBranchPoint} repeatedly, and
we will always use the map $f^{\circ n}$. The global arms at $z_1$
containing $0$ and $c_1$ are different because
$z_1 \in (0,c_1)$.  If the orbit is {\tame}, then the local arm $L_0$
cannot be mapped to itself; since $G_0 \owns 0$,
$L_0$ must map to $L_1$. There is a unique local
arm at $z_1$ which maps to the local arm towards $0$. Let
$G_{q-1}$ be the corresponding global arm; it may or may not map
onto $0$ under $f^{\circ k}$ for $k \leq n$. All the other
global arms are mapped onto their images homeomorphically. They can
be labelled so that $G_i$ maps to $G_{i+1}$ for $i=1,2,\ldots, q-2$.
This settles the {\tame} case.

In the {\evil} case, the local arm $L_0$ is fixed, and the other local arms
are permuted transitively. Let $L_{q-1}$ be the arm for which
$f^{\circ n}(L_{q-1})$ points to the critical value. Then all other global
arms map homeomorphically and can be labelled $G_1, G_2,\dots, G_{q-2}$
so that $G_i$ maps homeomorphically into $G_{i+1}$ for $i=1,2,\ldots, q-2$.

If $f^{\circ k}(G_{q-1}) \not\owns 0$ for all $k \leq n$,
then the entire cycle $G_1,\ldots,G_{q-1}$ of global arms
would map homeomorphically onto their images, and all their
endpoints would have identical itineraries with $z_1$. This
contradicts the expansivity condition for Hubbard trees.
\qed

\begin{coro}{Itinerary of Characteristic Point}
\label{CorCharactItin} \lineclear
In the Hubbard tree for the $\*$-periodic kneading sequence $\nu$,
fix a periodic point $z$ whose orbit does not contain the critical
point. Let $m$ be the period of $z$; if $z$ is not a branch point,
suppose that the itinerary of $z$ also has period $m$. Then if the
first $m-1$ entries in the itinerary of $z$ are the same as those
in $\nu$, the point $z$ is characteristic.

There is a converse if $z$ is a branch point: if $z$ is
characteristic, then the first $m$ entries in its itinerary are the
same as in $\nu$.
\end{coro}
\proof
If $z$ is not characteristic, then by
Lemma~\ref{LemCharacteristic}, the arc $[z,c_1]$ contains the
characteristic point of the orbit of $z$; call it $z'$. The
itineraries of $z$ and $z'$ differ at least once within the period
(or the period of the itinerary would divide the period of $z$; for
branch points, this would violate the expansivity condition, and
otherwise this is part of our assumption).
If the itinerary of $z$ coincides with $\nu$ for at least $m-1$ entries,
then the same must be true for $z'\in[z,c_1]$
(it is easy to check that for any Hubbard tree, the set of points sharing the same
$m-1$ entries in their itineraries is connected).  
Since the number of symbols $\0$ must be the same in the itineraries of $z$ and $z'$,
then $z$ and $z'$ must have identical itineraries, and this is a contradiction.

Conversely, if $z$ is the characteristic point of a branch orbit,
then by Lemma~\ref{LemHomeo}, $[z,c_1]$ maps homeomorphically onto
its image under $f^{\circ m}$ without hitting $0$, and the first
$m$ entries in the itineraries of $z$ and $c_1$ coincide.
\qed

The following result allows to distinguish {\tame} and {\evil} branch
points just by their itineraries.

\begin{prop}{Type of Characteristic Point}
\label{PropBranchType} \lineclear
Let $z_1$ be a characteristic periodic point. Let $\tau$ be the
itinerary of $z_1$ and let $n$ be the exact period of $z_1$. Then:
\begin{itemize}
\item
if $n$ occurs in the internal address of $\tau$, then the first
return map of $z_1$ sends the local arm towards $0$ to the local
arm toward $c_1$, and it permutes all local arms at $z_1$
transitively;
\item
if $n$ does not occur in the internal address of $\tau$, then the
first return map of $z_1$ fixes the local arm towards $0$ and
permutes all other local arms at $z_1$ transitively.
\end{itemize}
In particular, a characteristic periodic branch point of period $n$
is {\evil} if and only if the internal address of its itinerary does
not contain $n$.
\end{prop}
\proof
The idea of the proof is to construct certain precritical points $\z'_{k_j}\in[z_1,0]$ so that $[z_1,\z'_{k_j}]$ contains no precritical points $\z'$ with $\step(\z')\le\step(\z'_{k_j})$. Using these points, the mapping properties of the local arm at $z_1$ towards $0$ can be investigated. We also need a sequence of auxiliary points $w_i$ which are among the two preimages of $z_1$.

Let $\z'_1=0$ and $k_0=1$ and let $w_1$ be the preimage of $z_1$ that is contained in $T_{\1}$ and let $k_1\geq 2$ be maximal such that $f^{\circ(k_1-1)}|_{[z_1,w_1]}$ is homeomorphic. If $k_1<\infty$, then there exists a unique point $\z'_{k_1} \in (z_1,w_1)$ such that $f^{\circ k_1-1}(\z'_{k_1}) = 0$. All points on $[z_1, \z'_{k_1})$ have itineraries which coincide for at least $k_1-1$ entries.
If $k_1<n$ then the interval $[w_2,f^{\circ(k_1-1)}(z_1)]$ is non-degenerate
and contained in $f^{\circ(k_1-1)}((\z'_{k_1},z_1])$,
where $w_2$ denotes the preimage of $z_1$ that is not separated from $f^{\circ(k_1-1)}(z_1)$ by $0$. Let $y_{k_1}\in(\z'_{k_1},z_1)$ be such that $f^{\circ(k_1-1)}(y_{k_1})=w_2$.
Next, let $k_2 > k_1$ be maximal such that $f^{\circ(k_2-1)}|_{[z_1,y_{k_1}]}$ is homeomorphic.
If $k_2<\infty$, then there exists a point $\z'_{k_2} \in (z_1,\z'_{k_1})$ such that $f^{\circ k_2-1}(\z'_{k_2}) = 0$, and the points on $[z_1, \z'_{k_2})$ have the same itineraries for at least $k_2-1$ entries. If $k_2<n$ then, as above, the interval $[w_3,f^{\circ(k_2-1)}(z_1)]$ is non-degenerate (where again $w_3$ is an appropriate preimage of $z_1$) and there is a $y_{k_2}\in(\z'_{k_2},z_1)$ such that $f^{\circ(k_2-1)}(y_{k_2})=w_3$. Continue this way while $k_j<n$.

Note that the $\z'_{k_j}$ are among the precritical points on $[z_1,w_1]$ closest to $z_1$ (in the sense that for each $\z'_{k_j}$, there is no $\z'\in(z_1,\z'_{k_j})$ with $\step(\z')\le\step(\z'_{k_j})$; compare also Definition~\ref{DefPrecritical}), but the $\z'_{k_j}$ are {\em not} all precritical points closest to $z_1$; in terms of the cutting time algorithm, the difference can be described as follows: starting with $[z_1,w_1]$, we iterate this arc forward until the image contains $0$; when it does after some number $k_j$ of iterations, we cut at $0$ and keep only the closure of the part containing $f^{\circ k_j}(z_1)$ at the end. The usual cutting time algorithm would continue with the entire image arc after cutting, but we cut additionally at the point $w_{j+1}\in f^{-1}(z_1)$.

The point of this construction is the following:
let $\rho_{\tau}$ be the $\rho$-function with respect to $\tau$, i.e., $\rho_{\tau}(j) := \min\{ i > j \ : \ \tau_i \neq \tau_{i-j}\}$.
Then $k_1 = \rho_{\tau}(1)$ and $k_{j+1} =\rho_{\tau}(k_j)$ (if $k\ne 1$, then the exact number of iterations that the arc $[z_1,z_k]$ can be iterated forward homeomorphically is $\rho_\tau(k)-1$ times). Therefore we have constructed a sequence $\z'_{k_j}$ of precritical points on $[z_1,0]$ so that
(for entries less than $n$) $k_0,k_1,\dots=\orb_{\rho_\tau}(1)$, which is the internal address associated to $\tau$.

Recall that $n$ is the exact period of $z_1$. If $n$ belongs to the internal address, then there exists $\z'_n \in [z_1,0]$ and
$f^{\circ n}$ maps $[z_1, \z'_n]$ homeomorphically onto $[z_1,c_1]$. Therefore $f^{\circ n}$ sends the local arm towards
$0$ to the local arm towards $c_1$. By Lemma~\ref{LemHomeo}, $f^{\circ n}$ permutes all arms at $z_1$ transitively.

On the other hand, assume that $n$ does not belong to the internal address. Let $m$ be the last entry in the internal address before $n$. Then $f^{\circ (m-1)}$ maps $[z_1,\z'_m]$ homeomorphically onto $[z_{m},0]$, and the restriction to $[z_m,w_j]\subset[z_m,0]$ survives another $n-m$ iterations homeomorphically (maximality of $m$). There is a point $y_m\in[z_1,0]$ so that $f^{\circ (m-1)}([z_1,y_m])\to[z_m,w_j]$ is a homeomorphism, so $f^{\circ n}([z_1,y_m])\to[z_1,z_{n-m+1}]$ is also a homeomorphism. The local arm at $z_1$ to $0$ maps under $f^{\circ n}$ to a local arm at $z_1$ to $z_{n-m+1}$, and since $z_1$ is characteristic, this means that the local arm at $z_1$ to $0$ is fixed under the first return map.
The other local arms at $z_1$ are permuted transitively by Lemma~\ref{LemHomeo}.
\qed

\begin{defi}{Admissible Kneading Sequence and Internal Address}
\label{DefAdmissKneading} \lineclear
We call a $\*$-periodic kneading sequence and the corresponding
internal address {\em admissible}\index{admissible}
if the associated Hubbard tree contains no {\evil} orbit.
\end{defi}

This definition is motivated by the fact that a kneading sequence is admissible
if and only if it is realized by a quadratic polynomial; see below.
\hide{
, shown in
Proposition~\ref{PropAdmissAngle}, that a kneading sequence is
admissible if and only if there is an external angle which generates this
kneading sequence in the sense of Definition~\ref{DefKneading}.
Another equivalent statement is the following.
}

\begin{prop}{Embedding of Hubbard Tree}
\label{PropEmbedding} \lineclear
A Hubbard tree $(T,f)$ can be embedded into the plane so that $f$
respects the cyclic order of the local arms at all branch points if
and only if $(T,f)$ has no {\evil} orbits.
\end{prop}

\proof
If $(T,f)$ has an embedding into the plane so that $f$ respects the
cyclic order of local arms at all branch points, then clearly there can be no {\evil} orbit
(this uses the fact from Lemma~\ref{LemTree} that no periodic orbit of branch points contains a critical point).

Conversely, suppose that $(T,f)$ has no {\evil} orbits, so all local
arms at every periodic branch point are permuted transitively.
First we embed the arc $[0,c_1]$ into the plane, for example on a
straight line. Every cycle of branch points has at least its
characteristic point $p_1$ on the arc $[0,c_1]$, and it does not
contain the critical point. Suppose $p_1$ has $q$ arms. Take $s \in
\{ 1, \dots, q-1\}$ coprime to $q$
and embed the local arms at $p_1$ in such a way that
the return map $f^{\circ n}$ moves each arc over by $s$ arms
in counterclockwise direction.
This gives a single cycle for every $s < q$ coprime to $q$.
Furthermore, this can be done for all characteristic
branch points independently.

We say that two marked points $x,y$ are adjacent if $(x,y)$
contains no further marked point. If a branch point $x$ is already
embedded together with all its local arms, and $y$ is an adjacent marked
point on $T$ which is not yet embedded but $f(y)$ is, then draw a
line segment representing $[x,y]$ into the plane, starting at $x$
and disjoint from the tree drawn so far. This is possible uniquely
up to homotopy. Embed the local arms at $y$ so that $f\colon y\to
f(y)$ respects the cyclic order of the local arms at $y$; this is
possible because $y$ is not the critical point of $f$.

Applying the previous step finitely many times, the entire tree $T$
can be embedded. It remains to check that for every characteristic
branch point $p_1$ of period $m$, say, the map $f\colon p_1\to
f(p_1)=:p_2$ respects the cyclic order of the local arms. By
construction, the forward orbit of $p_2$ up to its characteristic
point $p_1$ is embedded before embedding $p_2$, and
$f^{\circ(m-1)}\colon p_2\to p_1$ respects the cyclic order of the
embedding. If the orbit of $p_1$ is {\tame}, the cyclic order induced
by $f\colon p_1\to p_2$ (from the abstract tree) is the same as the one induced by $f^{\circ
(m-1)}\colon p_2\to p_1$ used in the construction (already embedded in the plane), and the
embedding is indeed possible.
\qed

\remark
It is well known that once the embedding respects the cyclic order
of the local arms and their dynamics, then the map $f$ extends
continuously to a neighborhood of $T$ within the plane, and even to
a branched cover of the sphere with degree $2$.
See for example \cite{BFH}. This implies that the kneading
sequence of $T$ is generated by an external angle
and that $T$ occurs as the Hubbard tree of a quadratic polynomial
(compare \cite{BKS}).

It is not difficult to determine the number of different embeddings of $T$
into the plane (where we consider
two embeddings of a Hubbard tree into the plane as equal if
the cyclic order of all the arms at each branch point is the
same): if $q_1,q_2,\dots,q_k$ are the number of arms at the different
characteristic branch points and $\phi(q)$ is the Euler function counting
the positive integers in $\{1,2,\dots,q-1\}$ that are coprime to $q$, then
the number of different embeddings of $T$ respecting the dynamics is
$\prod_i \phi(q_i)$ \cite[Section~5]{BKS}. This also counts the number of
times $T$ is realized as the Hubbard tree of a postcritically finite quadratic
polynomial. If the critical orbit is periodic of period $n$,
then it turns out that $\prod_i \phi(q_i)<n$; see \cite[Section~16]{BKS}.

\hide{
: then Algorithm~\ref{AlgoTreeAngle} associates the same
external angles to them, and the extensions of the map to the plane
are equivalent in the sense of Thurston~\cite{DHT}. Different
embeddings yield different external angles and different Thurston
maps.
\begin{coro}{Number of Embeddings of Hubbard Tree}
\label{CorNumberEmbed} \lineclear
Let $(T,f)$ be a Hubbard tree without {\evil} branch points,
and let $q_1,q_2,\ldots$ be the number of arms of the different
characteristic branch points. Then there are $\prod \phi(q_i)$
different ways to embed $T$ into the plane such that $f$ extends to
a two-fold branched covering.
Here $\phi(q)$ is the Euler function counting the positive
integers $1 \leq i < q$ which are coprime to $q$.
\end{coro}
\proof
Follow the proof of Proposition~\ref{PropEmbedding}:
for each of the characteristic branch points $z_i$
there are $\ph(q_i)$ cyclic orders in which the $q_i$ arms can be
embedded. The choices can be made independently from all other
characteristic branch points. This gives the asserted number.
\qed
\remark
We show in Lemma~\ref{LemNumberEmbed2} that a
Hubbard tree with critical orbit of period $n$
has less than $n$ different embeddings into the plane which are
compatible with the dynamics.
}

\hide{
\begin{lemma}{Characteristic Points have Two Inverse Images}
\label{LemTwoInvImages} \lineclear
For every Hubbard tree $(T,f)$ with $\*$-periodic or preperiodic
kneading sequence, every characteristic periodic point $p\in
T \sm \{ c_1\}$ with itinerary $\tau(p) \neq \ovl{\1}$
has two different preimages
$\{p_0,-p_0\}=f^{-1}(p)\subset T$; these are separated by $0$,
i.e.,  $0\in(-p_0,p_0)$.
\end{lemma}
\proof
Since $f((\alpha,0))=(\alpha,c_1)\ni p$, there is always a preimage of $p$
in $(\alpha,0)\subset T_\1$.
Since $\tau(p)\not=\ovl{\1}$, there is a $p_i\in\orb(p)$ such that $p_i\in T_\0$.
Then $[p_i,0)\subset T_\0$ is mapped homeomorphically onto $[p_{i+1},c_1)$.
Since $p$ is characteristic, we have $p\in [p_{i+1},c_1)$,
and there is a point $p_0\in[p_i,0)\subset T_\0$ with $f(p_0)=p$.
\qed
}

\section{The Admissibility Condition}
\label{SecAdmissCondition}


\Intro{ In this section we derive the nature of periodic branch
points of the Hubbard tree from the kneading sequence
(Propositions~\ref{PropArmsNonOri} and \ref{PropArmsOri}). We
also prove a condition ({\em admissibility condition})
\index{admissibility condition} on the
kneading sequence which  decides whether there are {\evil} orbits:
Proposition~\ref{PropEvilCombinatorics} shows that an {\evil} orbit
violates this condition, and Proposition~\ref{PropArmsNonOri} shows
that a violated condition leads to an {\evil} orbit within the Hubbard
tree. Since a Hubbard tree can be embedded in the plane whenever
there is no {\evil} orbit (Proposition~\ref{PropEmbedding}), we obtain
a complete classification of admissible kneading sequences
(Theorem~\ref{ThmEvilAdmiss}).
}

\begin{defi}{The Admissibility Condition}
\label{DefAdmissCond} \lineclear
A kneading sequence $\nu\in\Syms$ {\em fails the admissibility
condition for period $m$}\index{admissibility condition}
if the following three conditions hold:
\begin{enumerate}
\item
the internal address of $\nu$ does not contain $m$;
\item
if $k<m$ divides $m$, then $\rho(k)\leq m$;
\item
$\rho(m)<\infty$ and if $r\in\{1,\ldots,m\}$ is congruent to $\rho(m)$ modulo $m$, then
$\orb_\rho(r)$ contains $m$.
\end{enumerate}
A kneading sequence {\em fails the admissibility condition} if it does so for some $m\ge 1$.

\noindent
An internal address {\em fails the admissibility condition} if its
associated kneading sequence does.
\end{defi}

\hide{
If a kneading sequence fails the admissibility condition for period $m$, then it follows that $\rho(m)$ is in the associated internal address (Lemma~\ref{LemCombiEndpoint}).
}

This definition applies to all sequences in $\Syms$, i.e., all sequences in $\{\0,\1\}^\Nplus$ and all $\*$-periodic sequences, provided they start with $\1$. However, in this section and the next we will only consider $\*$-periodic and preperiodic kneading sequences because these are the ones for which
we have Hubbard trees. The main result in this section is that this
condition precisely describes admissible kneading
sequences in the sense of Definition~\ref{DefAdmissKneading} (those
for which the Hubbard tree has no {\evil} orbits):

\begin{theo}{{\Evil} Orbits and Admissibility Condition}
\label{ThmEvilAdmiss} \lineclear
A Hubbard tree contains an {\evil} orbit of exact period $m$ if and only if
its kneading sequence fails the Admissibility Condition
\ref{DefAdmissCond} for period $m$.

Equivalently, a Hubbard tree can be embedded into the plane so that the dynamics
respects the embedding if and only if the associated kneading sequence does not fail
the Admissibility Condition \ref{DefAdmissCond} for any period.
\end{theo}
The proof of the first claim will be given in Propositions~\ref{PropEvilCombinatorics}
and \ref{PropArmsNonOri}, and the second is equivalent by Proposition~\ref{PropEmbedding}.

\hide{
While this theorem only deals with the postcritically finite case, we will show
in Corollary~\ref{CorAdmissAngle} that
a $\*$-periodic or non-periodic kneading sequence is admissible
in the sense of Definition~\ref{DefAdmissCond} if and only if it is generated
from an external angle by the algorithm in Definition~\ref{DefKneading}.
Note that the fact whether or not a sequence fails the admissibility condition for period $m$ is determined by its first $\rho(m)$ entries.
}

\begin{example}{Non-Admissible Kneading Sequences}
\label{ExNonAdmiss} \lineclear
The internal address $1 \IntAdr 2 \IntAdr 4 \IntAdr
5 \IntAdr 6$ with kneading sequence $\ovl{\1\0\1\,\1\0\*}$
(or any address that starts with  $1 \IntAdr 2 \IntAdr
4 \IntAdr 5 \IntAdr 6 \IntAdr$) fails the admissibility condition
for $m=3$, and the Hubbard tree indeed has a periodic branch point
of period $3$ that does not permute its arms transitively, as can
be verified in Figure~\ref{FigEvil}.
This is the simplest and best known example of a non-admissible
Hubbard tree; see \cite{IntAdr,Ke2,Pe}.

More generally, let $\nu = \ovl{\nu_1 \dots \nu_{m-1}\*}$
be any $\*$-periodic kneading sequence
of period $m$
so that there is no $k$ dividing $m$ with $\rho(k)=m$.
\hide{(for
kneading sequences in the Mandelbrot set, this means that $\nu$ is
not a bifurcation from a sequence of period $k$; see \ref{where?}).
}
This clearly implies
$\rho(k)<m$ for all $k$ dividing $m$.
Let $\nu_m \in \{ \0, \1\}$ be such that $m$ does not occur in the
internal address of $\ovl{\nu_1 \dots \nu_m}$.
Then for any $s\geq 2$, every sequence starting with
\[
\ovl{\underbrace{\nu_1 \dots \nu_m \dots \nu_1 \dots \nu_m}_{
\mbox{ $s-1$ times}} \nu_1 \dots \nu_{m-1}\nu'_m}
\]
(with $\nu'_m \neq \nu_m$) fails the admissibility condition for $m$.
The example $1 \IntAdr 2
\IntAdr 4 \IntAdr 5 \IntAdr 6$ above with kneading sequence
$\ovl{\1\0\1\1\0\*}$ has been constructed in this way, starting
from $\ovl{\1\0\*}$.

It is shown in \cite{AlexETDS} that every non-admissible kneading sequence
is related to such an example: the kneading sequences as constructed in this
example are exactly those where within the tree of admissible kneading sequences,
subtrees of non-admissible sequences branch off (compare also \cite[Section~6]{BKS}).
These sequences correspond exactly to primitive hyperbolic components of the
Mandelbrot set.

\hide{
We show after
Proposition~\ref{PropBranchNonAdmiss} that every non-admissible
kneading sequence is related to such an example.
}
\end{example}

\begin{figure}[htbp]
\begin{center}
\rule{22pt}{0pt}
\begin{minipage}{108mm}
\unitlength=9mm
\begin{picture}(12,3.2) \let\ts\textstyle
\put(2,2){\line(1,0){7}} \put(2,2){\circle{0.2}}
\put(1,1){\line(1,1){1}}
\put(1,3){\line(1,-1){1}}
\put(1.5,2.5){\line(1,1){0.5}}
\put(3,2){\line(0,-1){1}} \put(3,2){\circle{0.2}}
\put(4,2){\line(0,1){1}} \put(4,2){\circle{0.2}}
\put(5,2){\line(0,-1){1}}
\put(8,2){\line(0,1){1}} \put(8,2){\circle{0.2}}
\put(9,2){\line(1,1){1}} \put(9,2){\circle{0.2}}
\put(9,2){\line(1,-1){1}}
\put(9.5,1.5){\line(-1,-1){0.5}}
\put(6.2,2){\circle*{0.1}} \put(5.5,1.5){$c_0 = c_{11}$}
\put(1,1){\circle*{0.1}} \put(0.5,0.8){$c_1$}
\put(10,3){\circle*{0.1}} \put(10.2,2.8){$c_2$}
\put(4,3){\circle*{0.1}} \put(4.2,2.8){$c_3$}
\put(3,1){\circle*{0.1}} \put(2.5,0.8){$c_4$}
\put(8,3){\circle*{0.1}} \put(8.2,2.8){$c_5$}
\put(2,3){\circle*{0.1}} \put(2.2,2.8){$c_6$}
\put(9,1){\circle*{0.1}} \put(8.5,0.8){$c_7$}
\put(5,1){\circle*{0.1}} \put(4.5,0.8){$c_8$}
\put(1,3){\circle*{0.1}} \put(0.5,2.8){$c_9$}
\put(10,1){\circle*{0.1}} \put(10.2,0.8){$c_{10}$}
\end{picture}
\end{minipage}
\vskip-20pt
\end{center}
\caption{The Hubbard tree for $1 \IntAdr 2 \IntAdr 4 \IntAdr 5
\IntAdr 11$ is admissible. There is a {\tame} periodic orbit of branch points of
period $5$ (indicated by $\circ$'s). The other branch points are preperiodic.
}
\label{FigNotEvil}
\end{figure}
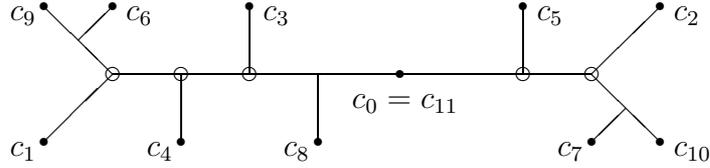

While $1 \IntAdr 2 \IntAdr 4 \IntAdr 5 \IntAdr 6$ is not
admissible, the internal address  $1 \IntAdr 2 \IntAdr 4 \IntAdr 5
\IntAdr 11$ is admissible; its Hubbard tree is shown in
Figure~\ref{FigNotEvil}.
This shows that the Translation Principle\index{Translation Principle}
from \cite[Conjecture~8.7]{IntAdr} does not hold: the address $1 \IntAdr
2 \IntAdr 4 \IntAdr 5 \IntAdr 11$ is realized in the $\frac13$ and
$\frac23$-sublimbs\index{sublimb}
of the real period $5$ component $1 \IntAdr 2
\IntAdr 4 \IntAdr 5$ of the Mandelbrot set. The Translation
Principle would predict that $1 \IntAdr 2 \IntAdr 4 \IntAdr 5
\IntAdr 6$ should exist within the $\frac12$-sublimb, but no such
hyperbolic component\index{hyperbolic component}
exists (the same counterexample was found
independently by V.~Kauko~\cite{VirpiPaper}).

\remark
The three conditions in the admissibility condition are independent:
here are examples of kneading sequences where exactly two of the
three conditions are satisfied.
\begin{itemize}
\item
$\nu = \ovl{101\*}$ ($1\IntAdr 2\IntAdr 4$), $m = 2$:
condition~1 is violated; $\nu$ is admissible.
\item
$\nu = \ovl{111\*}$ ($1\IntAdr 4$), $m = 2$: condition~2 is
violated; $\nu$ is admissible.
\item
$\nu = \ovl{101\*}$ ($1\IntAdr 2\IntAdr 4$), $m = 3$: condition~3
is violated; $\nu$ is admissible.
\end{itemize}

These conditions can be interpreted as follows: the first condition
picks a candidate period for an {\evil} orbit, taking into account that
a branch point is always {\tame} when its period occurs on the
internal address (Proposition~\ref{PropEvilCombinatorics}); the
second condition assures that the period $m$ of the {\evil} orbit is
the exact period, and the third condition makes the periodic orbit
{\evil} by assuring that the first return map of the characteristic
point maps a different local arm than the one pointing to $0$ onto
the local arm to the critical value.

\begin{lemma}{Bound on Failing the Admissibility Condition}
\label{LemBoundAdmissibility} \lineclear
If a $\*$-periodic kneading sequence of period $n$ fails the admissibility condition
for period $m$, then $m<n$.
\end{lemma}
\proof
Since $n$ occurs in the internal address, we may suppose $m\neq n$. If $m>n$, then $\rho(m)<m+n$ (because one of the entries between $m$ and $m+n$ is a $\*$), hence $r<n$ and $\orb_\rho(r)$ terminates at $n$, so $m\notin\orb_\rho(r)$.
\qed

A different way to interpret Lemma~\ref{LemBoundAdmissibility} is to say that a $\*$-periodic kneading sequence fails the admissibility condition for period $m$ if and only if the associated Hubbard tree has an evil branch point of period $m$ (Theorem~\ref{ThmEvilAdmiss}), and the period of a branch point is bounded by the period of the kneading sequence (Lemma~\ref{LemHomeo}).

\medskip

One of the main tools are closest precritical points.

\begin{defi}{Precritical Points}
\label{DefPrecritical} \lineclear
A point $x \in T$ is called {\em precritical} if $f^{\circ
k}(x)=c_1$ for some $k\geq 1$; the least such index $k$ is called
$\step(x)$\index{$\step$}.
The point $x$ is called a {\em closest precritical
point}\index{closest precritical point}
and denoted $\z_k$ if $f^{\circ j}([c_1,x])\not\ni c_1$ for
all $j\in \{1, \dots, k-1\}$.
\end{defi}

The critical point is always $\z_1$; if the critical point is
periodic of period $n$, then $\z_n=c_1$ and there is no closest
precritical point $x$ with $\step(x)>n$. Closest precritical
points are those which are ``visible from $c_1$'' in the sense of
\cite[Section 8]{IntAdr}: the idea is that a precritical point $\z$ blocks
the view of all $\z'$ behind $\z$ with $\step(\z')\ge\step(\z)$
(figuratively speaking, $\z$ is so big that the smaller point $\z'$
cannot be seen if it is behind $\z$).
We say that $\z$ is the {\em earliest} precritical point
on an arc $(x,y)$ (or $[x,y]$ etc.) if it is the one with the lowest $\step$.

\begin{lemma}{Closest Precritical Points Unique}
\label{LemClosestUnique} \lineclear
A Hubbard tree contains at most one closest precritical point
$\z_k$ for every index $k$.
\end{lemma}

\proof
If for some $k$, there are two closest precritical points $\z_k$
and $\z'_k$, then $f^{\circ (k-1)}$ maps $[\z_k,\z'_k]$
homeomorphically onto its image, but both endpoints map to the
critical point $0$. This is a contradiction.
\qed

\begin{lemma}{Elementary Properties of $\rho$}
\label{LemRho} \lineclear
If $\z_k\neq c_1$, then the earliest
precritical point on $(\z_k,c_1]$ is $\z_{\rho(k)}$.
For $k\geq 1$, the earliest precritical point on $[c_{1+k},c_1)$ is
$\z_{\rho(k)-k}$.

If $\z_k\neq c_1$, then
$[c_1,\z_k]$ contains those and only those closest precritical
points $\z_m$ for which $m\in\orb_\rho(k)$. In particular,
$\z_m\in[0,c_1]$ if and only if $m$ belongs to the internal address.
\end{lemma}

\proof
The first two statements follow immediately from the definition of
$\rho$, using the idea of {\em cutting times}, namely that if we look
at the largest neighborhood of $c_1$ in $T$ on which $f^{\circ k}$ is
monotone, we have cut this neighborhood at a closest precritical
point $\z_n$ whenever $n = k-1$.
Note that the arc $[c_1,\z_k)$ can be iterated
homeomorphically for at least $k$ iterations,
and $f^{\circ k}([c_1,\z_k)) = [c_{1+k}, c_1)$.
The first time that $f^{\circ(m-1)}([c_1,\z_k))$ hits $0$
is for $m=\rho(k)$ by definition,
and the earliest precritical point on $[c_1,\z_k)$
takes exactly $\rho(k)$ steps to map to $c_1$. The claim now
follows by induction. The statement about the internal address
follows because $\z_1 = 0$.
\qed

\begin{lemma}{Images of Closest Precritical Points}
\label{LemCPPAgain}\lineclear
If $k < k' \leq \rho(k)$, then $f^{\circ k}(\z_{k'})$
is the closest precritical point $\z_{k'-k}$.
\end{lemma}
\proof
Let $x := f^{\circ k}(\z_{k'})$.
Then the arc $[c_1,\z_{k'}]$ maps under $f^{\circ k}$
homeomorphically onto $[c_{k+1},x]$, and there is no precritical
point $\z\in(c_{k+1},x)$ with $\step(\z)\leq k'-k$.
If $\rho(k) > k'$ then by Lemma~\ref{LemRho} there is no such precritical point
$\z\in[c_1,c_{k+1}]$ either, and hence none on $(x,c_1]$. Since $\step(x)=k'-k$, the point $x$
is indeed the closest precritical point $\z_{k'-k}$.
Finally, if $\rho(k) = k'$, then $\z_{\rho(k)-k} = \z_{k'-k}$ is the
earliest precritical point on $[c_1, c_{1+k}]$ (Lemma~\ref{LemRho}).
But since $x$ also has $\step(x) = k'-k$ and $[x,\z_{k'-k}]$ contains no
point of lower $\step$, we have $x=\z_{k'-k}$.
\qed

\begin{lemma}{Precritical Points Near Periodic Points}
\label{LemPrecritNearPeriodic} \lineclear
Let $z_1$ be a characteristic periodic point of period $m$ such that
$f^{\circ m}$ maps $[z_1, c_1]$ homeomorphically onto its image.
Assume that $\nu$ is not $\*$-periodic of period less than $m$.
If $z_1$ has exactly two local arms, assume also that the first
return map of $z_1$ interchanges them.
Then
\begin{enumerate}
\item
the closest precritical point $\z_m$ exists in the Hubbard tree,
$z_1 \in [\z_m, c_1]$ and $\z_{\rho(m)} \in [c_1, z_1]$;
\item
if $\z$ is a precritical point closest to $z_1$ with $\step(\z) < m$
in the same global arm of $z_1$ as $\z_m$, then
$\z_m \in [z_1, \z]$;
\item if $z_1$ is a {\tame} branch point, then $\z_m \in [0,z_1]$ and $m$
occurs in the internal address;
\item if $z_1$ is an {\evil} branch point, then $\z_m \in G_{q-1}$
(where global arms are labelled as in Lemma~\ref{LemHomeo})
and $m$ does not occur in the internal address.
\end{enumerate}
\end{lemma}

\proof
(1) First we prove the existence of $\z_m$ in $T$.
Let $G_0, G_1, \dots, G_{q-1}$ be the global arms of $z_1$
with $0 \in G_0$ and $c_1 \in G_1$.
(Note that $q=2$ is possible.)
Let $L_0, \dots, L_{q-1}$ be the corresponding local arms.
Let $j$ be such that $f^{\circ m}(L_j) = L_1$.

If $j=0$, then $0=\z_1\in G_j$. If $j\neq 0$, then $q\geq 3$ by
assumption, so $j=q-1$ by Lemma~\ref{LemHomeo} and
there is an $i < m$ so that $f^{\circ i}(G_j)$ contains $0$.
Therefore, in both cases there
exists a unique $\z_k\in G_j$ with $k\leq m$ maximal, and it
satisfies $z_1\in(\z_k,c_1)$. We want to show that $k=m$.

If $k<m$, then $f^{\circ k}$ maps $(z_1,\z_k)$ homeomorphically
onto $(z_{k+1},c_1)\ni z_1$. By maximality of $k$, the restriction
of $f^{\circ m}$ to $(z_1,\z_k)$ is a homeomorphism with image
$(z_1,c_{m-k+1})\subset G_1$, and it must contain
$f^{\circ(m-k)}(z_1)=z_{1+m-k}$ in contradiction to the fact that
$z_1$ is characteristic. Hence $k=m$, $\z_m$ exists and $z_1 \in
[c_1,\z_m]$.

Clearly $f^{\circ m}$ maps $[z_1, \z_m]$ homeomorphically onto
$[z_1,c_1]$. By Lemma~\ref{LemRho}, $\z_{\rho(m)} \in [c_1, \z_m]$.
If $\z_{\rho(m)} \in [z_1, \z_m]$, then $f^{\circ m}(\z_{\rho(m)})
\in [c_1, z_1] \subset [c_1, \z_{\rho(m)}]$, but then
$\z_{\rho(m)}$ would not be a closest precritical point.
Hence $\z_{\rho(m)} \in [z_1, c_1]$.

(2) For the second statement, let $k:= \step(\z) < m$. We may suppose that $(z_1,\z)$
contains no precritical point $\z''$ with $\step(\z'')<m$ (otherwise replace $\z$ by $\z''$).
Clearly $\z \notin [z_1, \z_m]$. Assume by contradiction that
$\z_m\notin[z_1,\z]$ so that $[z_1,\z,\z_m]$ is a non-degenerate
triod. Since both $[z_1,\z]$ and $[z_1,\z_m]$ map homeomorphically
under $f^{\circ m}$, the same is true for the triod $[z_1,\z,\z_m]$.

Under $f^{\circ k}$, the triod $[z_1,\z,\z_m]$ maps homeomorphically
onto the triod $[z_{k+1},c_1, \z']$ with
$\step(\z')=m-k$, and $z_{k+1} \neq z_1$. Then
$z_1\in(z_{k+1},c_1)$, so the arc $(z_{k+1},\z')$ contains either
$z_1$ or a point at which the path to $z_1$ branches off.
Under $f^{\circ(m-k)}$, the triod $[z_{k+1},c_1,\z']\ni z_1$ maps homeomorphically
onto $[z_1,c_{m-k+1},c_1]\ni z_{m-k+1}$.
Therefore $(z_1,c_1)$ contains
either the point $z_{m-k+1}$ or a branch point from which the path to
$z_{m-k+1}$ branches off. Both are in contradiction to the characteristic
property of $z_1$.

(3) If $z_1$ is {\tame}, then $j=0$, so $\z_m \in G_0$.
By the previous statement, $\z_m \in [z_1,0]$, and
Lemma~\ref{LemRho} implies that $m$ belongs to the internal address.

(4) Finally, if $z_1$ is {\evil}, then $\z_m \in G_{q-1}$ and $m$ does not
occur in the  internal address by Lemma~\ref{LemRho}.
\qed

The following lemma is rather trivial, but helpful to refer to in longer arguments.
\begin{lemma}{Translation Property of $\rho$}
\label{LemRhoTranslate} \lineclear
If $\rho(m)>km$ for $k\ge 2$, then $\rho(km)=\rho(m)$.
\end{lemma}
\proof
Let $\nu$ be a kneading sequence associated to $\rho$. Then $\rho(m)>km$ says that the first $m$ entries in $\nu$ repeat at least $k$ times, and $\rho(m)$ finds the first position where this pattern is broken. By definition, $\rho(km)$ does the same, omitting the first $k$ periods.
\qed

\begin{lemma}{Bound on Number of Arms}
\label{LemNumArmsIntAddr} \lineclear
Let $z_1 \in [c_1, \z_m]$ be a characteristic point of period $m$
with $q$ arms.
Assume that $\nu$ is not $\*$-periodic of period less than $m$.
If $z_1$ has exactly two local arms, assume also that the first
return map of $z_1$ interchanges them.
If $z_1$ is {\evil}, then $(q-2)m<\rho(m) \leq (q-1)m$;
if not, then $(q-2)m<\rho(m) \leq qm$.
\end{lemma}
\proof
Let $G_0 \owns 0$, $G_1 \owns c_1, \dots, G_{q-1}$ be the global arms at
$z_1$. By Lemma~\ref{LemPrecritNearPeriodic}, $\z_{\rho(m)} \in
[z_1, c_1]$. The lower bound for $\rho(m)$ follows from
Lemma~\ref{LemHomeo}.

\begin{figure}[htb]
\begin{center}
\begin{minipage}{90mm}
\unitlength=9mm
\begin{picture}(10,3.2) \let\ts\textstyle
\put(0,2.5){\line(1,0){10}}
\put(0,2.5){\circle*{0.1}} \put(-0.1,2.7){$c_1$}
\put(2,2.5){\circle*{0.1}} \put(1.9,2.75){$\z_{\rho(m)}$}
\put(4,2.5){\circle*{0.1}} \put(3.9,2.7){$z_1$}
\put(8,2.5){\circle*{0.1}} \put(7.9,2.7){$0$}
\put(1,3){$G_1$}\put(6,3){$G_0$}\put(6.5,0){$T'$}

\put(4,2.5){\line(0,-1){1}}
\put(4,1.5){\circle*{0.1}} \put(4.2,1.7){$y$}

\put(4,1.5){\line(-3,-1){2}}
\put(3,1.15){\circle*{0.1}} \put(3,0.7){$\z_m$}

\put(4,1.5){\line(3,-1){2}}
\put(6,0.8){\circle*{0.1}} \put(6.1,1.05){$y'$}

\put(6,0.8){\line(1,0){3}}
\put(7,0.8){\circle*{0.1}} \put(7,1.05){$\z_r$}
\put(9,0.8){\circle*{0.1}} \put(8.6,1.05){$c_{1+(q-2)m}$}
\end{picture}
\end{minipage}
\vskip-20pt
\end{center}
\caption{Subtree with an {\evil} branch point $z_1$ of a Hubbard tree.}
\label{FigPointy}
\medskip
\end{figure}
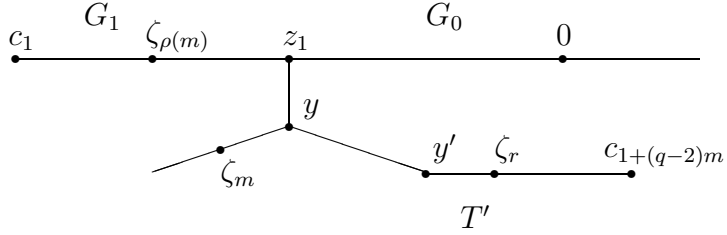

First assume that $z_1$ is {\evil}, so $q \geq 3$.
Lemma~\ref{LemRhoTranslate} implies $\rho((q-2)m)=\rho(m)$.
Assume by contradiction that $\rho(m) > (q-1)m$.
Then $r := \rho(m) - (q-2)m > m$. By
Lemma~\ref{LemCPPAgain}, $\z_r = f^{\circ ((q-2)m)}(\z_{\rho(m)})$
is a closest precritical point. It belongs to the same arm
$G_{q-1}$ as $\z_m$, but $\z_m \notin [z_1, \z_r]$.
As $\z_{\rho(m)}$ is the earliest precritical point on $[c_1, \z_m)$,
we cannot have $\z_r \in [z_1, \z_m]$ either.
Therefore $[z_1, \z_r, \z_m]$ is a non-degenerate triod within
$\ovl G_{q-1}$; let $y\in G_{q-1}$ be the branch point, see
Figure~\ref{FigPointy}.
Obviously
$f^{\circ m}(y) \in [z_1, c_1]$ and since $f^{\circ((q-2)m)}$ maps
$G_1$ homeomorphically into $G_{q-1}$, we find
$y' := f^{\circ ((q-1)m)}(y) \in [z_1, c_{1+(q-2)m}]$,
see Figure~\ref{FigPointy}.
If $y' \in [z_1, y]$, then $f^{\circ ((q-1)m}$ maps $[z_1, y]$
homeomorphically into itself. This contradicts expansivity of the tree.
Therefore $y' \in (y, c_{1+(q-2)m}]$. Let $T'$ be the component of
$T \setminus \{ y \}$ containing $c_{1+(q-2)m}$. Since
$\z_{\rho(m)} \in [z_1,c_1]$ and $f^{\circ (q-2)m}(\z_{\rho(m)}) = \z_r$,
$T'$ contains $\z_r$ but not $\z_m$.
Now $f^{\circ ((q-1)m)}$ maps $T'$  homeomorphically into itself
(otherwise, there would be an earliest precritical point $\z\in T'$ with
$\step(\z)<m$, but then $\z_m\in[z_1,\z]$ by
Lemma~\ref{LemPrecritNearPeriodic}~(2)). Again, expansivity of the
tree is violated. Thus indeed $\rho(m)\le (q-1)m$.

\begin{figure}[htb]
\begin{center}
\begin{minipage}{90mm}
\unitlength=9mm
\begin{picture}(10,6) \let\ts\textstyle
\put(0,2.5){\line(1,0){10}}
\put(0,2.5){\circle*{0.1}} \put(-0.1,2.7){$c_1$}
\put(2,2.5){\circle*{0.1}} \put(1.9,2.7){$\z_{\rho(m)}$}
\put(4,2.5){\circle*{0.1}} \put(3.9,2.7){$z_1$}
\put(5,2.5){\circle*{0.1}} \put(5.1,2.7){$y$}
\put(7,2.5){\circle*{0.1}} \put(6.9,2.7){$\z_m$}
\put(9,2.5){\circle*{0.1}} \put(8.9,2.7){$0$}

\put(1,3){$G_1$}\put(6,3){$G_0$}\put(3.9,0.7){$G_{q-1}$}
\put(4,2.5){\line(0,-1){1}}
\put(4,1.5){\circle*{0.1}}
\put(4,1.5){\line(-3,-1){2}}
\put(3,1.15){\circle*{0.1}} \put(3,0.7){$\z_k$}

\put(4,1.5){\line(3,-1){2}}

\put(6,0.8){\circle*{0.1}} \put(6.1,1){$c_{1+(q-2)m}$}

\put(5,2.5){\line(0,1){2.5}}
\put(5,3.5){\circle*{0.1}}  \put(5.1,3.4){$y'$}
\put(5,4.3){\circle*{0.1}}  \put(5.1,4.2){$\z_r$}
\put(5,5){\circle*{0.1}} \put(5.1,4.9){$c_{1+(q-1)m}$}

\put(5,3.5){\line(-1,1){1}}
\put(4,4.5){\circle*{0.1}} \put(3,4.7){$c_{1+m-k}$}
\end{picture}
\end{minipage}
\vskip-20pt
\end{center}
\caption{Subtree with a {\tame} branch point $z_1$ of a Hubbard tree.}
\label{FigPointy2}
\medskip
\end{figure}
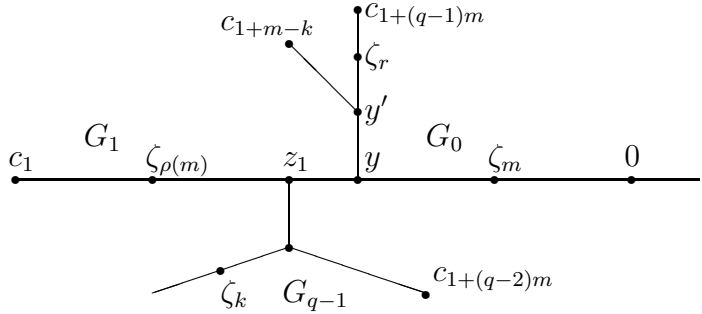

Now assume that $z_1$ is not {\evil} (and maybe not even a branch point).
Assume by contradiction that $\rho(m) > qm$.
We repeat the above argument with $r := \rho(m) - (q-1)m > m$,
conclude that $\rho((q-1)m)=\rho(m)$
and find the closest precritical point $\z_r \in G_0$, so $\z_r$ is in the same
global arm at $z_1$ as $\z_m$.
As before, $[z_1, \z_r, \z_m]$ is a non-degenerate triod with branch
point $y$ and $y' := f^{\circ qm}(y)$ lies on
$[y, c_{1+(q-1)m}]$.
Let $T'$ be the component of $T \setminus \{ y \}$ containing
$c_{1+(q-1)m}$.
We claim that $f^{\circ qm}$ is homeomorphic on $T'$.

It follows as above, using Lemma~\ref{LemPrecritNearPeriodic}, that
$f^{\circ m}$ maps $T'$ homeomorphically into $G_1$, and
$f^{\circ (q-2)m}$ maps $G_1$ homeomorphically into $G_{q-1}$.
Let $T'' = f^{\circ (q-1)m}(T')$ and
assume by contradiction that $f^{\circ m}$ is not homeomorphic on $T''$.
Then $T''$ contains a closest precritical point $\z_k$ for
some $k < m$. Take $k<m$ maximal.
Then $f^{\circ m}$ is homeomorphic on $[z_1,\z_k]$, and since
$f^{\circ k}([z_1,\z_k])=[z_{k+1},c_1]\ni z_1$, it follows that
$f^{\circ m}([z_1,\z_k])=[z_1,c_{1+m-k}]$ contains $z_{1+m-k}$. But
since $\z_k\in T''$, hence $f^{\circ(q-1)m}(y)\in [z_1,\z_k]$, we also
have $y'\in [z_1,c_{1+m-k}]$. As a result, $[z_1,c_{1+m-k}]\subset [z_1,y]\cup T'$.

If $z_{1+m-k}\in[z_1,y]$, then $f^{\circ m}$ maps $[z_1,z_{1+m-k}]$ homeomorphically
onto its image. This is a contradiction: both endpoints are fixed, but the image must be
in $G_1$.
Therefore, $z_{1+m-k}\in T'$. By
Lemma~\ref{LemPrecritNearPeriodic}~(2) again, there can be no
precritical point with $\step$ less than $m$ on $[z_1,z_{1+m-k}]$,
and we get the same contradiction.

We can conclude as above that $f^{\circ qm}$ maps $T'$
homeomorphically into itself as claimed.
But this is a contradiction
to expansivity of the tree.
\qed

\begin{prop}{{\Evil} Orbit Fails Admissibility Condition}
\label{PropEvilCombinatorics} \lineclear
If a Hubbard tree has an {\evil} orbit of exact period $m$
and $\nu$ is not $\*$-periodic of period less than $m$, then the
kneading sequence fails the admissibility condition for period $m$.
\end{prop}

\proof
Let $z_1$ be the characteristic point of the {\evil} orbit of period
$m$ and let $G_0, \dots , G_{q-1}$ be the global arms labelled as in
Lemma~\ref{LemHomeo}. The corresponding local arms will be labelled
$L_0,\ldots,L_{q-1}$.

We know from Lemma~\ref{LemPrecritNearPeriodic} that $m$ is not in
the internal address, so the first part of the admissibility
condition is already taken care of.
Since $[z_1,c_1]$ maps homeomorphically for $m$ steps, the first
$m$ entries in the itineraries of $z_1$ and $c_1$ coincide,
and $e(z_1) = \ovl{\nu_1 \dots \nu_m}$.
Let $k<m$ be a divisor of $m$. Suppose by contradiction that $\rho(k) > m$.
Then $e(z_1)$ has period $k$.

Therefore the period of $z_1$ is a multiple, and if it is a proper
multiple, then $z_1$ and $f^{\circ k}(z_1)$ are two periodic points
with the same itinerary.  This contradicts expansivity of the Hubbard tree,
so the period of $z_1$ must be $k$ as well. This however
contradicts the assumption, settling the second condition
of Definition~\ref{DefAdmissCond}

Let $r:=\rho(m)-(q-2)m$. By Lemma~\ref{LemNumArmsIntAddr},
$0<r\leq m$. By Lemma~\ref{LemCPPAgain},
$f^{\circ (q-2)m}(\z_{\rho(m)}) = \z_{\rho(m)-(q-2)m} =\z_r
\in G_{q-1}$.

By Lemma~\ref{LemPrecritNearPeriodic}~(4) and (2), we have $\z_m\in
G_{q-1}$ and then $\z_m \in [z_1, \z_r]$. Now Lemma~\ref{LemRho}
shows that $m \in \orb_{\rho}(r)$. Hence $\nu$ fails the
admissibility condition for period $m$.
\qed

In Propositions~\ref{PropArmsNonOri} and \ref{PropArmsOri}, we will
determine the exact number of arms at all branch points, and determine
from the internal address which branch points a Hubbard tree has.

\begin{prop}{Number of Arms at {\Evil} Branch Points}
\label{PropArmsNonOri}\lineclear
Suppose a kneading sequence $\nu$ fails the Admissibility
Condition~\ref{DefAdmissCond} for period $m$, and that
$\nu$ is not $\*$-periodic of period less than $m$. Then the Hubbard tree
for $\nu$ contains an {\evil} branch point of exact period $m$; the number
of its arms is
$q:=\lfloor \rho(m)/m\rfloor +2\geq 3$.
\end{prop}
\proof
Write $\rho(m)=(q-2)m+r$ for $r\in\{1,2,\ldots,m\}$ and $q\ge 3$.
Then $\rho((q-2)m)=\rho(m)$ by Lemma~\ref{LemRhoTranslate} and
the earliest precritical point on
$[c_{1+(q-2)m},c_1)$ is $\z_r$ by Lemma~\ref{LemRho}.
Since $\nu$ fails the admissibility condition for $m$, this
implies in particular $m\in\orb_\rho(r)$, hence by Lemma~\ref{LemRho}
$\z_m\in[\z_r,c_1]\subset[c_{1+(q-2)m},c_1]$.

Consider the connected hull
\[
H:=[c_1,c_{1+m},c_{1+2m},\ldots,c_{1+(q-3)m},\z_m] \,\,.
\]
Since $\rho(km)=\rho(m)>(q-2)m$ for $k=2,3,\dots,q-3$ by Lemma~\ref{LemRhoTranslate},
the map $f^{\circ m}$ sends the arc $[c_1,c_{1+km}]$
homeomorphically onto its image, and the same is obviously true for $[c_1,\z_m]$.
We thus get a homeomorphism $f^{\circ m}\colon H\to H'$ with
\[
H'=[c_{1+m},c_{1+2m},c_{1+3m},\ldots,c_{1+(q-2)m},c_1] \,\,.
\]
Since $\z_m\in [c_{1+(q-2)m},c_1]$, we have $H\subset H'\subset H\cup[\z_m,c_{1+(q-2)m}]$.
Moreover, $\z_r\in[c_{1+km},c_{1+(q-2)m}]$ for $k=0,1,\dots,q-3$:
the first difference between the itineraries of $c_{1+(q-2)m}$ and $c_1$
occurs at position $r$, while $c_1$ and $c_{1+km}$ have at least $m$ identical entries.
Since $\z_m\in[\z_r,c_1]$, it follows similarly that
$\z_m\in[\z_r,c_{1+km}]\subset[c_{1+(q-2)m},c_{1+km}]$ for $k\le q-3$.
Therefore, $H'\sm H=(\z_m,c_{1+(q-2)m}]$.

Among the endpoints defining $H$, only $c_{1+(q-3)m}$ maps outside
$H$ under $f^{\circ m}$, so $c_{1+(q-3)m}$ is an endpoint of $H$ and thus also of $H'$.
It follows that $c_{1+(q-4)m}$ is an endpoint of $H$ and thus also of $H'$ and so on,
so $c_1,\dots,c_{1+(q-2)m}$ are endpoints of $H$. Finally, also $\z_m$ is an endpoint
of $H$ (or $c_1$ would be an inner point of $H'$). As a result, $H$ and $H'$ have the
same branch points.

If $q=3$, then $H$ is simply an arc which is mapped in an orientation
reversing manner over itself, and hence contains a fixed point of
$f^{\circ m}$. Otherwise $H$ contains a branch point. Since $f^{\circ m}$
maps $H$ homeomorphically onto $H'\supset H$, it permutes
the branch points of $H$. By expansivity there can be at  most one
branch point, say $z_1$, which must be fixed under $f^{\circ m}$.
Since $f^{\circ m}\colon [z_1,c_1]\to[z_1,c_{1+m}]$ is a homeomorphism with
$[z_1,c_1]\cap[z_1,c_{1+m}]=\{z_1\}$, the arc $(z_1,c_1]$ cannot contain a
point on the orbit of $z_1$, so $z_1$ is characteristic.

If $z_1$ is a {\tame} branch point, then $\z_m\in[z_1,0]$ by
Lemma~\ref{LemPrecritNearPeriodic}, and $m$ occurs in the internal
address in contradiction to the failing admissibility condition.
If $z_1$ has exactly two arms, these are interchanged by
$f^{\circ m}$, and $\z_m\in G_0$, the global arm containing $0$.
By Lemma~\ref{LemPrecritNearPeriodic}~(2), $\z_m\in[0,z_1]$ and $m$
occurs in the internal address, again a contradiction. Hence
$z_1$ is an {\evil} branch point.

Now $H$ has exactly $q-1$ endpoints, and these are contained in
different global arms of $z_1$. The corresponding local arms are
permuted transitively by $f^{\circ m}$. Since $z_1$ is {\evil}, it has
exactly $q$ arms.
\qed

\noindent
This also concludes the proof of Theorem~\ref{ThmEvilAdmiss}.
\qed

\begin{defi}{Upper and Lower Kneading Sequences}
\label{DefUpperLower} \lineclear
If $\nu$ is a $\*$-periodic kneading sequence of exact period $n$,
we obtain two periodic kneading sequences $\nu_\0$ and $\nu_\1$ by
consistently replacing every $\*$ with $\0$ (respectively with
$\1$); both sequences are periodic with period $n$ or dividing
$n$, and exactly one of them contains the entry $n$ in its
internal address. The one which does is called the {\em upper
kneading sequence associated to $\nu$}\index{upper kneading sequence}
and denoted $\A(\nu)$, and
the other one is called the {\em lower kneading sequence associated
to $\nu$}\index{lower kneading sequence} and denoted $\Abar(\nu)$.
\end{defi}

\begin{lemma}{Itinerary Immediately Before $c_1$}
\label{LemImmediatelyBefore} \lineclear
When $x\to c_1$ in a Hubbard tree for the $\*$-periodic kneading sequence
$\nu$, the itinerary of $x$ converges (pointwise) to $\Abar(\nu)$.
\end{lemma}
\proof
Let $\tau$ be the limiting itinerary of $x$ as $x\to c_1$ and let $n$ be
the period of $\nu$. Then $\tau$ is clearly periodic with period
(dividing) $n$ and contains no $\*$, so $\tau\in\{\A(\nu),\Abar(\nu)\}$.
Let $m$ be the largest entry in the internal address of $\nu$ which is
less than $n$. Then there is a closest precritical point $\z_m\in[0,c_1)$  (Lemma~\ref{LemRho})
and $f^{\circ n}$ maps $[\z_m,c_1]$ homeomorphically onto its image.
Since $f^{\circ m}$ sends $(\z_m,c_1)\ni x$ onto
$(c_1,c_{1+m})\ni f^{\circ m}(x)$, which is sent by
$f^{\circ(n-m)}$ onto $(c_{1+n-m},c_{1+n})$, we get
$\tau_1\dots\tau_{n-m} = \nu_1\dots\nu_{n-m} = \tau_{m+1}\dots\tau_{n}$.
Hence $\rho_\tau(m)>n$; since $m$ occurs in the internal address of
$\tau$, the number $n$ does not.
\qed

\begin{prop}{Exact Period of Kneading Sequence}
\label{PropExactPeriod} \lineclear
For every $\*$-periodic kneading sequence of period $n$, the
associated upper kneading sequence $\A(\nu)$ has exact period $n$.
\end{prop}
\proof
Let $\tau:=\A(\nu)$ be the upper kneading sequence associated to
$\nu$ and suppose by contradiction that the exact period of $\tau$
is $m<n$. Then $\nu$ fails the admissibility condition for period
$m$: since $\rho_\tau(m)=\infty$ and $n$ is in the internal address
of $\tau$ by assumption, $m$ cannot occur on the internal address
of $\tau$ and hence neither on the internal address of $\nu$.
If $\rho_\tau(k)\geq m$ for a proper divisor $k$ of $m$, then the
exact period of $\tau$ would be less than $m$, a contradiction.
Hence $\rho_\nu(k)=\rho_\tau(k)<m$. The third part of the
admissibility condition is clear because $r=m$.

Thus by Theorem~\ref{ThmEvilAdmiss} the Hubbard tree for $\nu$, say
$(T,f)$, has an {\evil} orbit with period $m$. Let $z_1$ be its
characteristic point; it has itinerary $\tau$. Then $f^{\circ n}$
sends $[z_1,c_1]$ homeomorphically onto itself, so all points on
$[z_1,c_1)$ have itinerary $\tau$. By
Lemma~\ref{LemImmediatelyBefore}, it follows that $\tau$ is the
lower kneading sequence associated to $\nu$, a contradiction.
\qed

\hide{
A different proof for this result, based entirely on the combinatorics of
kneading sequences, is given in Lemma~\ref{LemExactPeriod}.
\remark
The exact period of $\Abar(\nu)$
can a proper divisor of $n$; see Lemma~\ref{LemPeriodicKneadingNoStar}.
}

\begin{lemma}{Characteristic Points and Upper Sequences}
\label{LemCharPtsUpperKS}\lineclear
Let $z$ be a characteristic point with itinerary $\tau$ and exact period $n$. Then exactly one of the following two cases holds:
\begin{enumerate}
    \item \begin{itemize}
        \item all local arms are permuted transitively, i.e., $z$ is tame,
        \item the internal address of $\tau$ contains the entry $n$,
        \item the exact period of $\tau$ equals $n$,
        \item $\tau=\A(\nu)$ for some $\star$-periodic kneading sequence $\nu$ of exact                 period $n$.
        \end{itemize}
    \item \begin{itemize}
        \item the local arm towards $0$ is fixed, all others are permuted transitively,
        \item the internal address of $\tau$ does not contain the entry $n$,
        \item if the exact period of $z$ and $\tau$ coincide then $\tau=\Abar(\nu)$ for some        $\star$-periodic kneading sequence $\nu$ of exact period $n$.
        \end{itemize}
\end{enumerate}
 For any $\star$-periodic sequence $\nu$, there is at most one tame periodic point in $T$ such that $\tau(p)=\A(\nu)$.
\end{lemma}

\proof
By Corollary~\ref{CorTwoKindsOrbits} either all local arms at $z$ are permuted transitively or the local arm pointing to $0$ is fixed and all others are permuted transitively.
By Proposition~\ref{PropBranchType}, $n$ is contained in the internal address of $\tau$ if and only of the local arm towards $0$ is not fixed.

Let $n'$ be the exact period of $\tau$. Then $n=kn'$ for some $k\geq 1$ and $\rho_\tau(n)=\infty$. Therefore, if $n$ is contained in the internal address
of $\tau$ then $k=1$ by the last assertion of Lemma~\ref{LemExactPeriod}.

The last remaining property of the first case follows immediately from $n=n'$, the definition of upper and lower kneading sequences and Proposition~\ref{PropExactPeriod}. Similarly, the third property in the second case follows from these results.

For the last statement, let us assume that there are two tame periodic points $p,q$ with itinerary $\A(\nu)$. Then they have both exact period $n'$ and $f^{\circ n'}([p,q])=[p,q]$. Thus not all local arms of $p,q$ are permuted transitively and neither $p$ nor $q$ is tame, a contradiction.
\qed

An immediate corollary of the preceding lemma is that if the period of $z$ and of $\tau$ coincide, then the type of $z$ (tame or not) is completely encoded in $\tau$.

In the second case however, if the exact period of $z$ and $\tau$ do not coincide, then $\tau$ may equal the upper or the lower kneading sequence of some $\star$-periodic kneading $\nu$ of exact period $n'$.

\begin{lemma}{Periodic Point behind Closest Precritical Point}
\label{LemPeriodicPrecrit} \lineclear
Let $\z\in[0,c_1)$ be a precritical point with $\step(\z) = m$
so that $f^{\circ m}:[\z,c_1] \to [c_1,c_{1+m}]$ is homeomorphic.
Then the arc $(\z,c_1)$ contains
a characteristic periodic point $z$ with exact period $m$. The
first return map of $z$ fixes no local arm at $z$.
\end{lemma}
\proof
First we show that $(c_1,\z)$ contains a periodic point of period $m$.
Assume by contradiction that this is not the case.
The construction in \cite{TreeExistence} does not only give the existence 
of the abstract Hubbard tree, but also the existence of extended trees 
that contain a finite number periodic orbits, see 
also \cite[Theorem 20.12]{BKS}.
Here we will include an $m$-periodic point $p$ with itinerary
$\tau = \ovl{\nu_1 \dots \nu_m}$, where $\nu$ is the
kneading sequence of the Hubbard tree.
If $c_1$ is periodic of period $n < m$, then we will
use the itinerary $\tilde \nu$ of a point $x$ very close to
$c_1$; so $\tau =  \ovl{\tilde \nu_1 \dots \tilde \nu_m}$.
By the choice of $\tau$, $\z$ does not separate $p$ from $c_1$.
\begin{figure}[htb]
\begin{center}
\begin{minipage}{90mm}
\unitlength=9mm
\begin{picture}(10,3.2) \let\ts\textstyle
\thicklines
\put(0,2.52){\line(1,0){6}} \put(0,2.5){\line(1,0){6}}
 \put(0,2.48){\line(1,0){6}}
\put(4,2.5){\line(0,-1){1.5}}\put(4.02,2.5){\line(0,-1){1.5}}
\put(3.98,2.5){\line(0,-1){1.5}}
\put(4,1){\circle*{0.1}} \put(4.2,0.8){$p$}
\thinlines
\put(0,2.5){\line(1,0){10}}
\put(0,2.5){\circle*{0.1}} \put(-0.1,2.7){$c_1$}
\put(2,2.5){\circle*{0.1}} \put(1.6,2.7){$f^{\circ m}(z)$}
\put(4,2.5){\circle*{0.1}} \put(3.9,2.7){$z$}
\put(6,2.5){\circle*{0.1}} \put(5.9,2.7){$\z$}
\put(8,2.5){\circle*{0.1}} \put(7.9,2.7){$0$}

\put(2,2.5){\line(0,-1){2.5}}
\put(2,0.8){\circle*{0.1}} \put(2.2,0.7){$f^{\circ 2m}(z)$}
\put(2,1.6){\circle*{0.1}} \put(2.2,1.5){$\z_r$}
\put(2,0){\circle*{0.1}} \put(2.2,0.0){$c_{1+m}$}
\put(1.4,1.3){$Y$}
\end{picture}
\end{minipage}
\vskip-20pt
\end{center}
\caption{Subtree $H_0 = [c_1,p,\z]$ (bold lines) with its image under
$f^{\circ m}$.
\label{FigPointy3}}
\medskip
\end{figure}
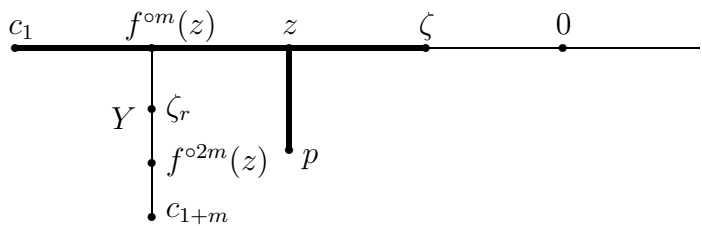
The triod $H_0 = [c_1,p,\z]$ maps under $f^{\circ m}$ homeomorphically
onto $[c_{1+m},p,c_1]$, see Figure~\ref{FigPointy3}.
If $H_0$ is degenerate, then it must necessarily have $c_1$ in the middle.
But then, $f^{\circ m}([c_1,p,\z]) = [c_{1+m},p,c_1]$ is degenerate
with $c_{1+m}$ in the middle: we have $\z\in[0,c_1]$, $c_1\in[\z,p]$ and
$c_{1+m}\in[c_1,p]$, hence $c_1\in[0,c_{1+m}]$ in contradiction to the fact
that $c_1$ is an endpoint of the Hubbard tree (we cannot have $c_{1+m}=c_1$
because then $\z=c_1$).

Hence there is a branch point, say $z$, in the interior of $H_0$.
Since $z\in(p,\z)$, we have $f^{\circ m}(z)\in(p,c_1)\subset [z,p)\cup[z,c_1)$.
The possibility $f^{\circ m}(z)=z$ contradicts our initial assumption.

If $f^{\circ m}(z) \in (z,p)$, then $f^{\circ m}$ maps $[z,p]$
homeomorphically into itself, so all points on $[z,p]$ have the same
itinerary. This contradicts either expansivity or finiteness of
the orbit of the branch point $z$.

Therefore, $f^{\circ m}(z) \in (c_1,z)$. In this case, $f^{\circ m}([c_1,z])$ branches
off from $[c_1,\z]$ at $f^{\circ m}(z)$; it belongs to an arm $Y$ at $f^{\circ m}(z)$,
and $f^{\circ 2m}(z) \in Y$. By expansivity, $f^{\circ m}$ cannot map $Y$
homeomorphically into itself, so there exists a closest precritical point $\z_k \in Y$
with $k < m$. By Lemma~\ref{LemRho}, $m \notin \orb_{\rho}(k)$. There is a unique $s \in
\orb_{\rho}(k)$ with $s < m < \rho(s)$. Then $f^{\circ m}$ maps the triod $[c_1,\z_s,z]$
homeomorphically onto the image triod $[c_{1+m}, c_{1+m-s}, f^{\circ m}(z)]$ with branch
point in $(f^{\circ m}(z),c_{1+m})\subset Y$. Therefore, $c_{1+m-s} \in Y$. Now let
$\z'$ be the earliest precritical point on $[c_1,c_{1+m-s}]$. By Lemma~\ref{LemRho},
$\step(\z')=\rho(m-s) - (m-s)$. The first assertion of Lemma~\ref{LemExactPeriod} \hide{
\footnote{This lemma is part of an analysis of the structure of $\rho$-orbits for
arbitrary $\0$-$\1$-sequences; it is independent of earlier results, so this forward
reference does not introduce a circular reference.} }
states that $m \in
\orb_{\rho}(\rho(m-s) - (m-s))$, so $\z_m \in [c_1,c_{1+m-s}]$. Therefore, $\z_m\neq\z$
(the points $\z$ and $\z_m$ are in different arms at $f^{\circ m}$), and this is a
contradiction.

We have now proved the existence of
a periodic point $z \in (c_1,\z)$ with itinerary $\tau$ and period $m$.
It is characteristic: if not, let $z_1\in(z,c_1)$ be the characteristic point; then
$f^{\circ m}(z_1,\z)=(z_1,c_1)$ and $z\in(z_1,\z)$, which is a contradiction.

Let $k|m$ be the exact period of $z$. By Lemma~\ref{LemHomeo}, $f^{\circ k}$
sends the local arm at $z$ to $0$ either to itself or to the local arm to $c_1$.
The first case is excluded by the fact that $f^{\circ m}\colon[z,\z]\to[z,c_1]$ is a homeomorphism.
In the second case, $f^{\circ k}\colon[z,\z]\to[z,f^{\circ k}(\z)]\subset[z,c_1]$ is a homeomorphism.
If $k<m$, then $f^{\circ k}(\z)\in(z,c_1)$ and $f^{\circ m}$ could not be a homeomorphism on $[\z,c_1]$. Hence $k=m$ is the exact period of $z$, and no local arm at $z$ is fixed by $f^{\circ m}$.
\qed

For any $m\geq 1$, let $r\in\{1,2,\ldots,m\}$ be
congruent to $\rho(m)$ modulo $m$, and define
\begin{equation}\label{numberarms}
q(m) := \left\{ \begin{array}{ll}
\frac{\rho(m)-r}{m} + 1 & \mbox{ if } m \in \orb_{\rho}(r)\,\,,\\
\frac{\rho(m)-r}{m} + 2 & \mbox{ if } m \notin \orb_{\rho}(r)\,\,.
\end{array} \right.
\end{equation}

\begin{prop}{Number of Arms at {\Tame} Branch Points}
\label{PropArmsOri}\lineclear
If $z_1$ is a {\tame} branch point of exact period $m$, then $m$ occurs in
the internal address, and the number of arms is $q(m)$.
Conversely, for any entry $m$ in the internal address with
$q(m)\geq 3$, there is a {\tame} branch point of exact period $m$ with
$q(m)$ arms (unless the critical orbit has period $m$).
\end{prop}

\proof
Let $z_1$ be the characteristic point of an orbit of {\tame}
branch points with exact period $m$,
and let $q'\geq 3$ be the number of arms at $z_1$.
By Lemma~\ref{LemHomeo}, the critical value cannot be periodic with period less than $m$.
By Lemma~\ref{LemPrecritNearPeriodic} (3) and (1),
$m$ occurs in the internal address, $\z_m\in(0,z_1)\subset G_0$, and $\z_{\rho(m)} \in [z_1,c_1]$.
Let $r' := \rho(m)-(q'-2)m$.
By Lemma~\ref{LemNumArmsIntAddr}, $0<r'\leq 2m$.
Therefore, by Lemma~\ref{LemRhoTranslate}, $\rho((q'-2)m)=\rho(m)$, so
by Lemma~\ref{LemCPPAgain}, $f^{\circ(q'-2)m}(\z_{\rho(m)}) = \z_{r'}$ is a
closest precritical point and by Lemma~\ref{LemHomeo},
$\z_{r'}=f^{\circ(q'-2)m}(\z_{\rho(m)})\in G_{q'-1}$.
Since $\z_m \in G_0$, we have
$\z_m \notin [c_1,\z_{r'}]$ and thus $m \notin \orb_{\rho}(r')$ (Lemma~\ref{LemRho}).

If $r' \leq m$, then $r = r'$ and we are in the case
$q(m) = \frac{\rho(m)-r}{m}+2 = \frac{\rho(m)-r'}{m}+2 = q'$.

If $r' > m$, then $r=r'-m$ and  $f^{\circ m}(\z_{r'})=\z_r$ by Lemma~\ref{LemCPPAgain}.
Then $f^{\circ m}$ maps $[z_1,\z_{r'}]$
homeomorphically onto $[z_1,\z_r]$, hence $\z_r\in G_0$.
By Lemma~\ref{LemPrecritNearPeriodic} (2) we find that
either $r=m$ or $r<m$ and $\z_m \in [z_1,\z_r]$, so in both cases $m \in \orb_{\rho}(r)$.
Therefore, $q(m) = \frac{\rho(m)-r}{m} + 1 = q'$. Again $q(m)=q'$.

\medskip

For the converse, let $m$ be an entry in the internal address.
By Lemma~\ref{LemRho}, the closest precritical point $\z_m$ exists on $[0,c_1]$.
By Lemma~\ref{LemPeriodicPrecrit}, $\z_m$ gives rise to a characteristic point
$z_1 \in [\z_m,c_1]$ of exact period $m$
and the first return map of $z_1$ fixes no local arm.
By Lemma~\ref{LemPrecritNearPeriodic} (1), $\z_{\rho(m)} \in [z_1,c_1]$.

If $z_1$ is a branch point, then
no local arm of $z_1$ is fixed by $f^{\circ m}$, so $z_1$ is {\tame}.
By the first assertion of the lemma, the number of arms is $q(m)$.

Finally, suppose that $z_1$ has only two arms $G_0 \owns 0$ and $G_1 \owns c_1$.
If the critical orbit is periodic and $m$ is an entry in the internal address, the period of the
critical orbit is at least $m$, and equality is excluded by hypothesis. By Lemma~\ref{LemNumArmsIntAddr}, we have $\rho(m)\le 2m$ and $r':=\rho(m)-m=r$. Then
$f^{\circ m}(\z_{\rho(m)}) = \z_{r'} \in G_0$.
Lemma~\ref{LemPrecritNearPeriodic} (2) then gives that $\z_m \in [\z_r,z_1]$
and hence $m \in \orb_{\rho}(r)$.
It follows that $q(m) = \frac{\rho(m)-r}{m} + 1 = 2$.
\qed

Together, Propositions~\ref{PropArmsNonOri} and
\ref{PropArmsOri} describe all branch points in all Hubbard trees.

\begin{coro}{Uniqueness of Hubbard Tree}
\label{CorUniqueness} \lineclear
If $(T,f)$ is a Hubbard tree with $\*$-periodic or preperiodic kneading sequence $\nu$, then $\nu$ alone determines $(T,f)$ uniquely up to equivalence.
\end{coro}
\proof
By Propositions~\ref{PropEvilCombinatorics} and \ref{PropArmsNonOri}, the tree $(T,f)$ has an evil periodic orbit of exact period $m$ if and only if $\nu$ fails the admissibility condition for period $m$; the number of arms is determined by Proposition~\ref{PropArmsNonOri}. By Proposition~\ref{PropArmsOri}, there is a branch point of period $m$ only if $m$ occurs in the internal address associated $\nu$; for every $m$ on this internal address, the quantity $q(m)$ from (\ref{numberarms}) determines whether or not there is a branch point and, if so the number of arms. Every branch point of any period $m$ has the property that the itinerary of the associated characteristic point coincides with $\nu$ for at least $m$ entries; this determines the itinerary of all points on the orbit of every branch point. Finally, every endpoint of $(T,f)$ is on the critical orbit by definition, so the itineraries of endpoints are shifts of $\nu$. 

If $(T',f')$ is another Hubbard tree with kneading sequence $\nu$, then we show that it is equivalent to $(T,f)$ in the sense as defined after Definition~\ref{DefHubbard}. Itineraries define a bijection between branch points of $(T,f)$ and $(T',f')$ and between postcritical points, and this bijection is respected by the dynamics. It thus suffices to prove that both trees have the same endpoints and their edges connect corresponding points. Recall that postcritical points and branch points are jointly known as marked points.

To see this, we use precritical points: by definition, these are points $\z\in T$ with $f^{\circ k}(\z)=c_0$ for some $k\ge 0$; in this case we write $\step(\z)=k$. Every such precritical point $\z$ has itinerary $\tau_1\tau_2\dots\tau_{k-1}\*\nu$ with $\tau_1,\dots,\tau_{k-1}\in\{\0,\1\}$. By induction on $k$, we show that such a point $\z$ with itinerary $\tau_1\tau_2\dots\tau_{k-1}\*\nu$ exists in $T$ if and only if it exists in $T'$, and if it does, corresponding marked points are in corresponding components of $T\sm\{\z\}$ resp.\ $T'\sm\{\z\}$. This is obvious for $k=0$ and $\z=c_0$. 

A precritical point $\z$ with $\step(\z)=k+1$ (described by the first $k$ entries $\tau_0,\dots,\tau_k$ of its itinerary) obviously exists if and only if there are two postcritical points $x,y\in T$ with $\z\in[x,y]$. This is equivalent to the existence of two points $x',y'\in T$ which are either postcritical points or precritical points with $\step(x)\le k$, $\step(y)\le k$ so that $\z\in[x',y']$ and so that $(x',y')$ does not contain a precritical point $\z'$ with $\step(\z')<k+1$. The latter condition can be checked using the itineraries of $x'$ and $y'$, and their existence is the same for $T$ and for $T'$ by inductive hypothesis. 

It now follows easily that $T$ and $T'$ have endpoints with identical itineraries, so they have a natural bijection between marked points. It also follows that every precritical point $\z$ disconnects $T$ and $T'$ into two parts so that corresponding parts contain marked points with identical itineraries, and this implies that the edges of $T$ and $T'$ connect corresponding points. This means by definition that $(T,f)$ and $(T',f')$ are equivalent as claimed.
\qed

\hide{
In \cite{IntAdr}, a special sufficient condition for admissibility of kneading sequences was derived and used to determine the Galois groups of periodic points: see Section~\ref{Sec:IntAddr}. We show that this condition follows easily from our general Admissibility Condition~\ref{DefAdmissCond}.
}

\end{document}